\documentclass[11pt]{article}
\usepackage{amssymb,amsmath,amsthm,amsfonts,mathtools}
\usepackage{mathrsfs}
\usepackage[bbgreekl]{mathbbol}
\usepackage[utf8x]{inputenc}
\usepackage{anyfontsize}
\usepackage[T1]{fontenc}

\usepackage{url}
\usepackage{float}

\usepackage{lmodern}
\usepackage{anyfontsize}
\usepackage{float} 
\usepackage[labelfont=bf]{caption}
\usepackage{placeins}
\usepackage{url}
\usepackage{bbding,bbold,bbm,empheq}
\usepackage[colorlinks=true, linkcolor=blue, citecolor=red]{hyperref}
\usepackage{cleveref}
\crefname{equation}{}{}
\crefname{thm}{}{}
\crefname{lem}{}{}
\crefname{defn}{}{}
\crefname{rk}{}{}
\crefname{prop}{}{}
\crefname{coro}{}{}
\crefname{section}{}{}
\crefname{subsection}{}{}

\usepackage{tikz, graphicx}
\usetikzlibrary{patterns}
\usepackage{caption, subcaption}
\usepackage{pdfpages}

\usepackage{titlesec}

\titleformat{\paragraph}[block]
  {\normalfont\normalsize\bfseries}
  {\theparagraph}{1em}{}

\titlespacing*{\paragraph}
  {0pt}{3.25ex plus 1ex minus .2ex}{1ex}
  
\usepackage{array}
\usepackage{placeins}
\usepackage{hyperref} 
\hypersetup{
    colorlinks=true,
    linkcolor=blue,      
    filecolor=black,
    urlcolor=black,
    citecolor=red,
    pdfpagemode=FullScreen,
}
\numberwithin{equation}{section}
\usepackage{dsfont}
\usepackage{fullpage}
\usepackage{float}

\usepackage{enumitem}
\setlist[itemize]{label=\textbullet}

\usepackage{authblk}
\usepackage[bottom]{footmisc}
\usepackage{amsthm}  
\makeatletter
\def\tagform@#1{\maketag@@@{(\textcolor{blue}{#1})}}
\makeatother
\usepackage{xcolor}
\usepackage{refcount}  

\renewcommand{\eqref}[1]{%
  (%
  \textcolor{blue}{\getrefnumber{#1}}%
  )%
}

\usepackage{color}
\usepackage{lmodern}

\usepackage{ucs}
\usepackage{tikz}
\usetikzlibrary{decorations.pathmorphing}
\theoremstyle{plain}
\newtheorem{thm}{Theorem}[section]
\newtheorem{prop}[thm]{Proposition}
\newtheorem{lem}[thm]{Lemma}

\theoremstyle{definition}
\newtheorem{defn}[thm]{Definition}

\newtheorem{ass}[thm]{Assumption}

\newcommand{\QED}{\mbox{}\hfill \raisebox{-0.2pt}{\rule{5.6pt}{6pt}\rule{0pt}{0pt}} \medskip\par}

\newenvironment{Proof}
{\noindent \abovedisplayskip = 0.5\abovedisplayskip
 \belowdisplayskip=\abovedisplayskip{\bfseries Proof. }}{\QED\medskip}

\title{\textbf{Supercritical  Geometric Optics 
\\ for  multidimensional scalar conservation laws}\\[0.5ex]%
       }
\author{
Christophe Cheverry
\thanks{Universit\'{e} de Rennes, IRMAR, CNRS, France, christophe.cheverry@univ-rennes.fr},
Rida Harb 
\thanks{Universit\'{e} C\^ote d'Azur, LJAD, CNRS, Nice, France, Rida.Harb@univ-cotedazur.fr},
St\'{e}phane Junca
\thanks{Universit\'{e} C\^ote d'Azur, LJAD, CNRS, Nice, France, Stephane.Junca@univ-cotedazur.fr}
}
\date{}

\begin{document}

\maketitle

\abstract{This article studies the propagation, near a smooth  non-constant ground state, of monophase
oscillatory wave trains with small amplitude and very high frequency for
multidimensional nonlinear scalar conservation laws. These oscillations lie
beyond the classical regime of uniformly bounded variation. Two scenarios may
occur. Under suitable compatibility conditions, the oscillations propagate.
Otherwise, they are damped through shock formation. In the baseline
supercritical regime, we give an explicit criterion distinguishing these two
behaviors and justify the corresponding asymptotic descriptions. At higher
frequencies, the propagation mechanism becomes more restrictive and is governed
by compatibility conditions involving a Lie algebra. Under a natural local
constant-rank assumption, we obtain a necessary and sufficient condition for
the existence of complete WKB expansions. We also prove cancellation under
first-order incompatibility. Perturbing a non-constant background state limits
the propagation of small scales compared with the constant-background case.
Our constructions contribute to the understanding of the fractional smoothing
effect conjectured by Lions, Perthame and Tadmor.}

\vspace{1em}
{\small
\medskip


\medskip

\noindent {\bf Data availability statement}: the data are available on request from any author.

\tableofcontents

\section{Introduction}\label{introduction}
The propagation of high-frequency waves in nonlinear hyperbolic equations raises a number of questions about the interaction between rapid oscillations, the nonlinear effects, and shock formation. In the multidimensional scalar setting, this issue is particularly interesting because highly oscillatory data may lie beyond the classical uniformly bounded $BV$ regime, while the entropy solution remains globally defined \cite{dafermos2016hyperbolic}. Our purpose is to understand under which conditions such oscillations can persist, and how they are instead damped by the nonlinear dynamics. More generally, nonlinear geometric optics provide a useful framework for quantifying the regularity of entropy solutions and for better understanding the compactness properties of the nonlinear semigroup generated by Kružkov entropy solutions. Precisely, given a smooth flux $F\in C^\infty(\mathbb{R},\mathbb{R}^d)$, we consider the  multidimensional scalar conservation law
\begin{equation}\label{eq:claw}
\partial_t u + \operatorname{div}_x F(u)=0,
\qquad u\in \mathbb{R},\quad x\in \mathbb{R}^d,\quad 1\leq d .
\end{equation}
Let $\underline{u}_0 \in C^\infty_b(\mathbb{R}^d,\mathbb{R})$ and 
$ \phi_0 \in C^\infty(\mathbb{R}^d,\mathbb{R}) $ be smooth bounded functions. Let $ \mathbb{T}:=\mathbb{R}/\mathbb{Z} $ be 
the torus. We select a profile  $U_0\in C^\infty_b(\mathbb{R}^d\times\mathbb{T},\mathbb{R})$ which is 
periodic with period one with respect to the second variable
$\theta\in\mathbb{T} $. 

\smallskip

\noindent We want to incorporate non-stationary phases $ \phi_0 $ and non-trivial profiles $ U_0 $, for which
\begin{equation}\label{nonstationary}
\nabla_x \phi_0 \not \equiv 0 , \qquad \partial_\theta U_0 \not \equiv 0 .
\end{equation}
Given $p\in\mathbb{N}^* $ and $ \varepsilon \in]0,1] $, at time $t=0$, we impose the initial condition
\begin{equation}\label{eq:initdata}
u^\varepsilon\big|_{t=0}
=u_0^\varepsilon(x)
:=\underline{u}_0(x)
+\varepsilon\,U_0\!\left(
x,\frac{\phi_0(x)}{\varepsilon^p}
\right).
\end{equation}
Our purpose is to describe the asymptotic behavior when $ \varepsilon $ goes to zero of the solutions $ u^\varepsilon $ to the Cauchy problem \eqref{eq:claw}--\eqref{eq:initdata}.

\smallskip
\noindent
The case $p=1$ corresponds to the classical framework of weakly nonlinear
geometric optics \cite{JMR}. Then, the first-order derivatives of the family $ (u_0^\varepsilon)_\varepsilon $ are  uniformly bounded, and these uniform
Lipschitz estimates allow the construction
of smooth solutions on a small time interval $[0,T]$ with $T>0$ independent of
$\varepsilon$. This holds true even for systems, when $ u \in \mathbb R^n $ with $ n \geq 1 $. The asymptotic description on $[0,T]$ of the corresponding solutions $ (u^\varepsilon)_\varepsilon $ has motivated a wide range of papers, including~\cite{DM,HMR,JMR}.

\smallskip
\noindent The situation $p>1$ is much more singular. This regime is usually referred to as \textit{supercritical}.  Under \eqref{nonstationary}, the lifespan $ T^\varepsilon >0 $ of {\it smooth} solutions $ u^\varepsilon $ generically shrinks to $ 0 $ when $ \varepsilon $ goes to $ 0 $, typically like $ \varepsilon^{p-1} $. This phenomenon already occurs when looking at the Burgers' equation, for $ d = 1 $ and $ F(u) = u^2/2 $. In the absence of specific assumptions, it applies along the same lines when $ n \geq 1 $ or $ d \geq 1 $. The construction of {\it weak} solutions for systems ($ n > 1$) is in general problematic\footnote{For instance, when $ d = 1 $, the existence theory for systems of two (or more) conservation laws generally requires BV-estimates. But, when $ p>1 $, the BV-norm of $ u_0^\varepsilon $ explodes when $ \varepsilon \rightarrow 0 $. As a consequence,the lifespan of these BV-solutions is expected to go to $ 0 $.}. However, when $ n = 1 $, this difficulty can be circumvented by implementing the notion of $ L^\infty $-entropy solution developed by Kruzkov~\cite{Kr}. Observe that
\[ M :=  
\sup_{0<\varepsilon\leq 1}\|u_0^\varepsilon\|_{L^\infty} < + \infty .
\]
Thus, provided shock formation is accounted for,
the Cauchy problem \eqref{eq:claw}-\eqref{eq:initdata} leads to global-in-time  $ L^\infty $-solutions. Still, to our knowledge, the asymptotic behavior of $ (u^\varepsilon)_\varepsilon $ has not yet been fully identified, even though some progress has been made in this direction \cite{CJ1,CJ2,CJR,J4}. Our aim is to address this gap. We work locally in space, near the origin $ x=0 $. Given $\rho>0$, we localize at time $ t= 0 $ the data to the ball of radius $\rho $, that is 
$ \Omega_0:=\{x\in\mathbb{R}^d:\ |x|\leq \rho\} $.
The speed of propagation of information is finite, limited by
\[ 
L:=\sup_{|u|\leq M}|a(u)|, \qquad a(u) := F'(u).
\]
Given $ T >0 $ with $ L \, T \leq \rho $, the domain of determinacy issued from $\Omega_0$ is defined by
\begin{equation}\label{Omega}
\Omega:=\{(t,x):\ t\in[0,T],\ |x|\leq \rho-Lt\} .
\end{equation}
For $t\in[0,T]$, we denote by $\Omega_t$ the section of $\Omega$ at time $t$, that is the ball 
\[
\Omega_t:=\{(t,x):\ |x|\leq \rho-Lt\} .
\]
The solution $\underline{u}$ to \eqref{eq:claw} with initial data $ \underline{u}_0 $ is locally smooth on $ \mathbb R_+ \times \mathbb R^d $, near $ (t,x) = (0,0) $. By restricting $ T >0 $ if necessary, we can always assume that $\underline{u}$ is smooth  on $\Omega$. From now on, we fix such a $ T $. Retain that $\underline{u} \in C^\infty_b (\Omega) $ is the smooth solution to
\begin{equation}\label{eq:transportform}
\partial_t \underline{u}
+a(\underline{u})\cdot\nabla_x\underline{u}=0,
\qquad
(t,x)\in\Omega,
\qquad
\underline{u}(0,x)=\underline{u}_0(x).
\end{equation}
The phase $\phi$ is transported according to the eikonal equation
\begin{equation}\label{eq:eikonal}
\partial_t\phi
+a(\underline{u})\cdot\nabla_x\phi=0,
\qquad
\phi(0,x)=\phi_0(x).
\end{equation}
Supercritical regimes appear when  $p>1$. Then, one cannot expect the oscillations to  propagate in every case. Indeed, such oscillations are no longer uniformly bounded in $BV$, so that two different scenarios are expected. In the first one, the oscillations persist and are propagated, but this requires strong compatibility conditions, which become increasingly restrictive when $p>2$. In the second one, the nonlinear dynamics of the Kružkov semigroup destroys the oscillatory structure, leading instead to a damping or averaging of the oscillations.

\smallskip
\noindent Some of these scenarios are stated rigorously in Section~\ref{sec:1+}, while Section~\ref{sec:2} provides a complete justification of the case $p=2$. Section~\ref{sec:3} is devoted to the more complex regime of higher frequencies $\varepsilon^{-p}$, with $p>2$, where a similar dichotomy persists. In this setting, the compatibility conditions become more complicated and more restrictive, and are expressed through the Lie algebra $\mathcal A$ generated by the vector fields $L_k:=a^{(k)}(\underline{u}_0)\cdot\nabla_x$, with $1 \leq k\leq p-1$. We identify a necessary and sufficient condition for the existence of a valid WKB geometric-optics expansion at any order, as stated in Theorem~\ref{thm:complete-WKB-criterion}. We also prove a cancellation result in Theorem~\ref{thm-cancel3}.

\smallskip
\noindent
Our results are related to the regularizing properties of the nonlinear semigroup associated with the scalar conservation law, as discussed in Appendix~\ref{app:regularity-compactness-semigroup}. Indeed, following the kinetic formulation
of Lions, Perthame and Tadmor~\cite{LPT}, several works have investigated fractional smoothing
effects, their optimality, and the limitations of $BV$ regularization in several space dimensions, 
see in particular~\cite{ACCM26,GGJJ,GJ}. From this perspective, the construction of highly oscillating solutions provides a natural test to evaluate :
\begin{enumerate}[label=\roman*)]
    \item The maximal regularity that can be expected from the nonlinear semigroup ;
    \item The mechanisms of cancellation, when they occur.
\end{enumerate}

\smallskip
\noindent
In connection with the research direction i), WKB expansions~\cite{CJ1} or highly singular data~\cite{CJ2,DW} have  already been used to exhibit
families of solutions involving critical fractional regularity. In particular, under a genuine nonlinearity assumption and near a
constant ground-state~\cite{CJR}, nontrivial oscillating ansatz with frequency $\varepsilon^{-p}$ have been
constructed as long as $p\leq d$, and this dimensional restriction turns
out to be optimal~\cite{J4}. One outcome of the present work is to show that the picture is deteriorating substantially when $\underline{u} $ is no longer constant\footnote{The study of what happens near a non-constant ground state $ \underline u $  is a way to investigate the regularity properties of scalar conservation laws in the presence of the variable (in $ t $ and $ x $) flux functions $ \tilde F(t,x,u) = F(\underline u +u)- F(\underline u) $.}. For a constant ground-state, the compatibility conditions
involve constant vector fields\footnote{It follows that the compatibility conditions reduce to simple algebraic orthogonality conditions.} which leave 
enough freedom to construct nontrivial phases in the range described above. 

\smallskip
\noindent
By contrast, for a
non-constant ground-state, the compatibility conditions become differential constraints coupling
the phase $\phi_0$ with the spatial variations of $\underline{u}_0$. Further compatibility conditions appear (when $ p > 2 $), and their geometric structure involves not
only the vector fields $ L_k $ but also the induced Lie algebra $ \mathcal A $ whose dimension may be strictly larger than $ p-1 $. Generically, the
additional directions produced by the Lie brackets force a compatible phase to be locally constant, obstructing the propagation of nontrivial oscillations. Thus, compared with the constant
ground-state setting, the variations of $ \underline u $ reveal supplementary  geometric obstructions for the propagation of supercritical oscillations\footnote{It follows that the limitation on $ p $ for the propagation of oscillations is in general strictly  less than $ d $ when $ d \geq 3 $.}.

\smallskip
\noindent
To investigate ii), the difficulty is to handle the formation of shocks, which do occur and are necessary for the cancellation process. The idea is to examine what happens  in a time boundary layer (of size $ \varepsilon $ when $ p= 2 $ which means to deal with $ \tau $ defined by $ \varepsilon \, \tau = t $), and to exploit \cite{Ch92,Ch97} in order to justify there a WKB expansion.
There is a hidden profile $ V(\tau,x,\theta) $ satisfying a Burgers law with respect to $ \theta $, developing discontinuities in $ \theta $, and decreasing in $ \tau $ according to Oleinik's estimates until it becomes negligible (for sufficiently large values of $ \tau $).

\section{Main results} \label{sec:1+}

We start with $p=2$ to explain main ideas and already obtain new results. Then, we present the outcomes for $p>2$, involving more complicated constraints related to a Lie algebra. 

\subsection{The basic supercritical regime}

The case $p=2$ is important in itself. It corresponds already to the propagation of oscillations that are not bounded in either $BV$ or $W^{1,1}$.
As will be seen, for $ p = 2 $, a complete discussion can be carried out. 
The case $p=2$ can serve to introduce the more complicated situation $p>2$.
Below, we present two main results, both obtained in Section~\ref{sec:2} and related to the supercritical regime $p=2$.
The first statement describes what happens when, at the initial time, the vector $a'(\underline{u}_0)$ is everywhere orthogonal to the gradient of $\phi_0$.

\begin{ass}[Compatibility condition between $ \underline{u}_0 $ and $ \phi_0 $] \label{ass-comp} The initial ground-state $\underline{u}_0$ and the initial phase $\phi_0$ are adjusted in such a way that
\begin{equation}\label{eq:compat-p2}
a'(\underline{u}_0)\cdot\nabla_x\phi_0\equiv 0
\qquad\text{on }\Omega_0 .
\end{equation}
\end{ass}

\noindent When $ d= 1 $ and the phase $ \phi_0 $ is nonstationary, this may happen only if the function $ a' (\underline{u}_0) $ vanishes, that is in the absence of nonlinearity. The situation is more favorable when $ d>1 $. Indeed, given any $ \underline{u}_0 $, whatever the choice of $ a $ with $ a' \bigl(\underline{u}_0(0) \bigr) \not = 0 $, it suffices to prescribe the trace of $ \phi_0 $ on an  hyperplane transversal to the direction $ a' \bigl(\underline{u}_0(0) \bigr) $ in order to recover (locally) a phase $ \phi_0 $ solving the transport equation \eqref{eq:compat-p2}.
We will prove in Lemma \ref{Lemma propagation} that the relation 
\eqref{eq:compat-p2} is propagated at the level of $ (\underline{u},\phi) $ on $ \Omega $. This property can be viewed as a geometric  transparency condition\footnote{The notion of transparency comes into play when the effects of the nonlinearity are weaker than expected.} which eliminates from the WKB calculus a problematic term of size one. Then, in factor of the power $ \varepsilon $, a formal computation reveals the modulation equation
\begin{equation}\label{Profile-intro}\partial_t U+\operatorname{div}_x \bigl(U\,a(\underline{u})\bigr)
+\frac{1}{6}\,
\bigl(a''(\underline{u})\cdot\nabla_x\phi\bigr) \,
\frac{\partial U^{3}}{\partial\theta}
=0,
\end{equation}
which must be associated with the initial data
\begin{equation}\label{Profile-introini}
U(0,x,\theta)=U_0(x,\theta).
\end{equation}
For $T>0$ sufficiently small, the Cauchy problem \eqref{Profile-intro}--\eqref{Profile-introini} has a smooth solution on $ \Omega \times \mathbb T $.

\begin{thm}[Propagation of oscillations uniformly bounded in $W^{1/2,1}$ but not in $ W^{s,1}$ for $ s >1/2$]
\label{propagation-thm-intro} Let $\phi$ be as in \cref{eq:eikonal} and $ U $ be as in 
\eqref{Profile-intro}--\eqref{Profile-introini}. Then, under Assumption \ref{ass-comp}, the lifespan $ T^\varepsilon $ of the  smooth solution $ u^{\varepsilon} $ to the initial value problem \eqref{eq:claw}-\eqref{eq:initdata} where $ p = 2 $ is uniformly bounded below by some $ T > 0 $. Moreover, on the domain $ \Omega $ given by \eqref{Omega},
in the sense of the sup-norm, we have
\begin{equation}\label{eq:expansion-osc}
u^{\varepsilon}(t,x)
=
\underline{u}(t,x)
+\varepsilon\,U\!\left(
t,x,\frac{\phi(t,x)}{\varepsilon^{2}}
\right)
+O(\varepsilon^{2}).
\end{equation}
\end{thm}

\smallskip
\noindent
Some comments are in order. This kind of propagation is always possible in dimension $d>1$, just by involving compatible data\footnote{Note that, when $ p=2 \leq d $, the Lie algebra generated by the  $L^k$ with $ 1 \leq k \leq 1 $ is of dimension $ 1 $, induced by $ L_1 $.
 The condition $ N_2=1 < d $ required in Theorem \ref{thm:complete-WKB-criterion} is systematically satisfied.}. 
  But, in dimension $d=1$, such propagation is possible  only if the flux is not uniformly convex.
Note that the case $p=2$ already contains several of the techniques used later for $p>2$. In particular, the cancellation phenomenon is proved first for $p=2$ and then extended to the case $p>2$.

\smallskip
\noindent
The converse scenario\footnote{What happens near a transition zone, that is when  $ a'(\underline{u}_0)\cdot\nabla_x\phi_{0 \mid x=0} = 0$ without having \eqref{eq:compat-p2} on $ \Omega_0 $ for some positive $ \rho $ small enough, will not be investigated.} arises when the orthogonality \eqref{eq:compat-p2} fails at  one point, $x=0$ for instance, and thereby (for $ \rho $ small enough) on all $ \Omega_0 $.

\begin{ass}[The genuinely nonlinear situation] \label{ass-fullyn} The initial ground-state $\underline{u}_0$ and the initial phase $\phi_0$ are adjusted in such a way that
\begin{equation}\label{violation-intro}
\inf_{\Omega_0} \ 
\bigl|a'(\underline{u}_0)\cdot\nabla_x\phi_0\bigr|>0.
\end{equation}
\end{ass}

\noindent Under \eqref{violation-intro}, instead of \eqref{Profile-intro}, we have to consider the (linear) transport equation
\begin{equation}\label{Profile-2-intro}
\partial_t\overline{U}
+\operatorname{div}_x \bigl(\overline{U}\,a(\underline{u})\bigr)
=0 .
\end{equation}
At time $ t = 0 $, this must be completed by 
the mean value with respect to $ \theta $ of $ U_0 $, that is
\begin{equation}\label{Profile-2-intromean}
\overline{U}(0,x)
=
\overline{U}_0(x)
:=
\int_0^1 U_0(x,\theta)\,d\theta.
\end{equation}
On the domain of determinacy $ \Omega $ yielding $ \underline{u} $, the Cauchy problem \eqref{Profile-2-intro}--\eqref{Profile-2-intromean} has a smooth solution $ \overline{U}(t,x) $ on $ \Omega  $.

\begin{thm}[Cancellation of oscillations]
\label{cancellation-thm-intro} Let $ \overline{U} $ be as in \eqref{Profile-2-intro}--\eqref{Profile-2-intromean}. Denote by $ u^{\varepsilon} $ the weak solution to \eqref{eq:claw}-\eqref{eq:initdata}  where $ p = 2 $. Then, under Assumption \ref{ass-fullyn}, for all $ t \in[0,T] $, in the sense of the $ L^1 $-norm on the ball $ \Omega_t $, we have
\begin{equation}\label{eq:expansion-avg}
u^\varepsilon(t,x)
=
\underline{u}(t,x)
+\varepsilon\,\overline{U}(t,x)
+o(\varepsilon) .
\end{equation}
\end{thm}

\smallskip
\noindent This theorem asserts that if the first compatibility condition is not satisfied, then the oscillations are canceled. Implicitly, Assumption~\ref{ass-fullyn} means that the vector field $L_1$ never vanishes. This is a genuinely nonlinear assumption on the flux along the ground state $\underline{u}_0$ and the phase $ \phi_0 $.  

\subsection{Higher frequencies}
We consider here the supercritical regime of oscillations with amplitude
$\varepsilon$ and frequency $\varepsilon^{-p}$, with $p>2$. We study
solutions of \cref{eq:claw} issued from initial data having the form 
\begin{equation}\label{eq:initdata-p>2}
u_0^\varepsilon(x)
=
\underline u_0(x)
+ \sum_{k=1}^K \varepsilon^k \,U_0^k \!\left(x,\frac{\phi_0(x)}{\varepsilon^p}\right)\, , \qquad K \in \mathbb N^* \, 
.
\end{equation}
In this section, our aim is to precisely characterize the conditions under which such oscillations persist over a time interval while maintaining a similar WKB expansion, as described in \eqref{eq:complete-WKB-expansion}. Before moving forward, let us introduce some key vocabulary. 

\begin{defn}[Non-trivial phase]
\label{def:nontrivial-phase}
A phase
$\phi\in C^\infty(\Omega;\mathbb R)$ is said to be
\emph{non-trivial} when
\begin{equation}\label{eq:nontrivial-phase-p}
\nabla_x\phi(t,x)\neq0,
\qquad
\forall (t,x)\in\Omega.
\end{equation}
\end{defn}

\noindent Oscillations cannot be detected when $ U(t,x,\cdot) $ is constant\footnote{Interesting situations appear only on condition that $ \partial_\theta U \not \equiv 0 $.}. This leads to the following notion.

\begin{defn}[Non-trivial profile]
\label{def:nontrivial-phase-profile-p}
A profile
$U\in C^\infty(\Omega\times\mathbb T;\mathbb R)$ is said to be
\emph{non-trivial} when
\begin{equation}\label{eq:nontrivial-profile-p}
\forall (t,x)\in\Omega,\qquad
\exists \, \theta\in\mathbb T \, ; \qquad 
\partial_\theta U(t,x,\theta)\neq0.
\end{equation}
\end{defn}

\noindent In what follows, a crucial role is played by the geometry of the following $ p-1 $ vector fields 
\begin{equation}\label{eq:time-dependent-generators}
\mathbb R^d \supset \Omega_t \ni x \longmapsto L_k(t,x)
:=
a^{(k)} \bigl( \underline u(t,x) \bigr)\cdot\nabla_x \, , \qquad 1 \leq k \leq p-1 .
\end{equation}

\begin{defn}[Compatible phase]
\label{def:compatible-phase}
A phase
$\phi\in C^\infty(\Omega;\mathbb R)$ is said to be
\emph{compatible} at time $ t $ (up to the order $ p $ for the background state $ \underline u $) when
\begin{equation}\label{imposeLk}
(L_k \, \phi) (t,x) = 0 \, , \qquad \forall \, x \in \Omega_t \, , \qquad 1 \leq k \leq p-1 .
\end{equation}
\end{defn}

\noindent In other words, at each point $(t,x)\in\Omega$, the vector fields $ L_k $ generate the generalized 
distribution
\begin{equation}\label{eq:distribution-Dp}
\mathcal D_p(t,x)
:=
\operatorname{span}
\bigl\{
L_1(t,x),\ldots,L_{p-1}(t,x)
\bigr\}
\subset T_x\mathbb R^d \, , 
\end{equation}
and the gradient of $ \phi $, that is $ \nabla_x \phi(t,x) $,  must be orthogonal to $ \mathcal D_p(t,x) $. It is therefore easy, by  composing vector fields, to deduce that $ X \phi \equiv 0 $ for all vector fields $ X $ belonging to the Lie closure of $ D_p(t,\cdot) $, namely
\begin{equation}\label{eq:compatibility-algebra-p}
\mathcal A_p(t,\cdot)
:=
\operatorname{span}
\left\{
X \, ; \, 
X\in
\operatorname{Lie}
\bigl(
L_1(t),\ldots,L_{p-1}(t)
\bigr)
\right\} \supset \mathcal D_p(t,x) \, .
\end{equation}
The distribution $ 
\mathcal A_p(t,\cdot) $ is involutive. In the real analytic case, by Sussmann-Nagano Theorem \cite{Nagano66}, it is integrable and generates the leaves of a foliation \cite{Fro} corresponding to 
the level surfaces of phases $ \phi $. But, in the actual $ C^\infty $-framework, to recover the property of integrability through  Frobenius' Theorem, some supplementary  condition is needed on the rank of $ \mathcal A_p(t,\cdot) $. Define
\begin{equation}\label{eq:def-Np}
\mathbb N \ni N_p(t,x)
:=
\dim\mathcal A_p(t,x) 
\leq d \, . 
\end{equation}
The function $ N_p (t,\cdot) : \Omega_t \rightarrow \mathbb N $ is always lower semi-continuous. It is required to be regular, and therefore locally constant. In particular, this must hold at time $ t = 0 $.

\begin{ass}[Locally constant rank of the Lie  algebra induced by $ \underline u $]
\label{ass:constant-rank-Np}
There exists $\rho>0$ sufficiently small such that the map $ x\mapsto N_p(0,x) $ 
is constant on $\Omega_0$. 
\end{ass}

\noindent Under Assumption \ref{ass:constant-rank-Np}, we just denote by $ N_p $ the (local near the origin) constant value of $N_p (0,x)$. Retain that
\begin{equation}\label{eq:constant-rank-Np}
N_p(0,x)=N_p,
\qquad \forall \, 
x\in\Omega_0.
\end{equation}
Note that this constant rank condition is imposed only at time $ t = 0 $, on $ \Omega_0 $ instead on all  $ \Omega_t $. This is because it is  preserved by the
evolution of the smooth ground state, see Lemma~\ref{lem:propagation-rank-Ap}.  As a consequence,
the same rank $N_p$ is valid throughout the domain $\Omega$. At this stage, it is important to pin down what happens when $ p = 2 \, $:
\begin{itemize}
\item When $ a' \circ \underline u \equiv 0$, the conservation law is linearly degenerate along $ \underline u $. We find $ N_2 = 0 $, and any smooth phase $ \phi $ is  compatible \cite{CGM03}.
\item When $ a' \circ \underline u (0,0) \not = 0 $, restricting $ \rho $ if necessary, we still have $ a' \circ \underline u (t,x) \not = 0 $ for all $ (t,x) \in \Omega $. There is only one non-zero vector field $ L_1(t,x) \not = 0 $, so that
\[ 
\mathcal D_2(t,x) = \mathcal A_2 (t,x) = \bigl \lbrace \lambda \,  L_1(t,x) \, ; \, \lambda \in \mathbb R  \bigr \rbrace \, , \qquad N_2 (t,x) = 1 \, ,  \qquad \forall \, (t,x) \in \Omega \, . \]
This means that $ N_2 = 1 $. As seen in Section \ref{sec:2}, as soon as $ 2 \leq d $, a non-trivial phase $ \phi $ (in fact a family of $ d-1 $ linearly independent phases $ \phi $) satisfying \eqref{imposeLk} do exist. 
\item The intermediate situation, that is when $ a' \circ \underline u $ is zero at $ (0,0) $ without being identically zero in a neighbourhood of the origin, is  not investigated in Section \ref{sec:2}. This corresponds in fact to a jump at $ x = 0 $ of $ N_2 (0,\cdot) $, breaking Assumption \ref{ass:constant-rank-Np}.
\end{itemize}

\noindent It is also necessary to grasp what happens when $ \underline u $ is constant. Then, $ \mathcal D_p = \mathcal A_p $ is constant, as well as $ N_p = \dim \mathcal D_p $. Assumption \ref{ass:constant-rank-Np} is automatically satisfied.

\smallskip

\noindent This short discussion is instructive because it indicates, firstly,  that the adequate geometrical criterion is based on $ \mathcal A_p $ instead of $ \mathcal D_p $ (though this distinction is not relevant when $ p= 2 $ or when $ \underline u $ is constant) and secondly, why Assumption \ref{ass:constant-rank-Np} is crucial (even when $ p= 2 $).

\smallskip

\noindent Observe that $ N_p (t,x) = d $ locally near the origin if and only if $ N_p (0,0) = d $. Thus, when $ N_p (0,0) = d $, 
Assumption \ref{ass:constant-rank-Np} is verified with $ N_p = d $. In PDEs, this  condition $ N_p = d $ is usually referred to as H\" ormander's  condition\footnote{In control theory, it is called the LARC (Lie Algebra Rank Condition) or Chow-Rashevskii condition \cite{Coron07}.}. 
It destroys any possibility of finding a non-constant phase $ \phi $ constrained by \eqref{imposeLk}. And therefore, it is a  geometrical obstruction to the persistence of supercritical oscillations.

\smallskip

\noindent When $ \underline u $ is constant, it is necessary and sufficient to impose $ \dim \mathcal D_p < d $ in order to exhibit a non-trivial compatible phase\footnote{When the vector fields $ L_k (t,0) $ are linearly independent, we obtain that $ \dim \mathcal D_p = p-1 $. In this setting, supercritical oscillations like in \eqref{eq:initdata-p>2} can propagate if and only if $ p \leq d $.}.
However, when $ 2 < p $ and $ \nabla_x \underline u \not \equiv 0 $, it is well known  that, generically, $ \mathcal D_p \subsetneqq \mathcal A_p $, and in general that $ N_p = d $ \cite{Coron07,Hormander1967}. This explains why the situation $ 2 < p $ together with the presence of variations of $ \underline u $ is much less favourable to the propagation of supercritical oscillations. 
We can now state the main propagation result of this section, which generalizes Theorem \ref{propagation-thm-intro}. It gives a
necessary and sufficient geometrical condition for the existence of
non-trivial complete supercritical WKB expansions at frequency
$\varepsilon^{-p}$.

\begin{thm}[Necessary and sufficient condition for complete supercritical WKB expansions]
\label{thm:complete-WKB-criterion}
Under Assumption~\ref{ass:constant-rank-Np} on $ \underline u $, there exist a radius
$\rho>0$, a time $T>0$, a non-trivial smooth compatible phase
$\phi\in C^\infty(\Omega;\mathbb R)$, a non-trivial principal profile
$U^1\in C^\infty(\Omega\times\mathbb T;\mathbb R)$ and, for every
$K\in\mathbb N^*$, profiles
$U^2,\ldots,U^K\in C^\infty(\Omega\times\mathbb T;\mathbb R)$
such that solutions to \eqref{eq:claw} exist with
\begin{equation}\label{eq:complete-WKB-expansion}
u^\varepsilon(t,x)
=
\underline u(t,x)
+
\sum_{k=1}^{K}
\varepsilon^k \, 
U^k\!\left(
t,x,\frac{\phi(t,x)}{\varepsilon^p}
\right)
+
\varepsilon^{K+1} \, r^\varepsilon(t,x),
\end{equation}
together with
\begin{equation}\label{eq:complete-WKB-remainder}
\|r^\varepsilon\|_{L^\infty(\Omega)}
\lesssim 1,
\qquad
\|\varepsilon^p \, \nabla_{t,x}r^\varepsilon\|_{L^\infty(\Omega)}
\lesssim 1,
\end{equation}
if and only if the Lie algebra $ \mathcal A_p $ is of dimension $ N_p<d $.
\end{thm}

\medskip
\noindent
Theorem~\ref{thm:complete-WKB-criterion} 
characterizes the existence of supercritical WKB expansions by a
purely geometrical condition on $ \underline u $. The rank $N_p$ measures the number of
spatial directions generated by the fields
$L_1,\ldots,L_{p-1}$ and all their successive Lie brackets, which are all involved during the construction process behind \eqref{eq:complete-WKB-expansion}. If $N_p<d$, at least one
transverse direction remains available, which makes it possible to
find a non-trivial phase annihilated by the Lie 
algebra. By contrast, if $N_p=d$, the Lie closure fills the
spatial tangent space. As a consequence, any phase annihilated by $\mathcal A_p$ must
satisfy $\nabla_x\phi=0$. Thus the condition $N_p<d$ expresses exactly
the geometric freedom needed for the persistence of non-trivial
supercritical oscillations.

\smallskip
\noindent
The propagation part of Theorem \ref{thm:complete-WKB-criterion} (when $ N_p<d $) involves compatible phases. We shall also consider the case of a phase which violates the first condition $ (L_1 \phi) \equiv 0 $ inside \eqref{imposeLk}. This means to impose  Assumption \ref{ass-fullyn} in order to extend Theorem \ref{cancellation-thm-intro} as indicated below.

\begin{thm}[Cancellation of oscillations]
\label{thm-cancel3}
Let $\overline U$ be as in \eqref{Profile-2-intro}-\eqref{Profile-2-intromean}. Denote by $ u^{\varepsilon} $ the weak solution to \eqref{eq:claw}-\eqref{eq:initdata}  where $ p \geq  2 $. Then, under
Assumption~\ref{ass-fullyn}, in the sense of the $L^1$-norm, we have
\begin{equation}\label{eq:expansion-avg-p>2}
u^\varepsilon(t,x)
=
\underline u(t,x)
+
\varepsilon\,\overline U(t,x)
+
o(\varepsilon).
\end{equation}
\end{thm}

\section{ Amplitude $\varepsilon$ with $\varepsilon^{-2}$ frequency}\label{sec:2}

\noindent
This section is devoted to the  supercritical regime obtained for the choice of $p=2$, corresponding
to initial data \cref{eq:initdata} which are bounded in $W^{1/2,1}$.
Under the compatibility condition \cref{eq:compat-p2}, the oscillations are
propagated, whereas a uniform violation of this condition leads to their
cancellation. These two situations are studied successively in
Subsections~\ref{Propagation} and~\ref{Cancellation}.
\par\medskip

\subsection{Propagation}\label{Propagation}

\medskip
\noindent
We establish here the propagation result stated in
Theorem~\ref{propagation-thm-intro}. The proof is organized in two steps.
In Subsubsection~\ref{propagation-compatibility}, we first show that the
compatibility condition imposed at the initial time is propagated by the
evolution. In Subsubsection~\ref{geometric-optics-reduction}, we then derive
the geometric-optics reduction, introduce the exact $\varepsilon$-dependent
profile, and conclude the proof by a classical stability argument.

\noindent
\subsubsection{Propagation of the compatibility condition}
\label{propagation-compatibility}

The compatibility condition \cref{eq:compat-p2} is imposed only on the
initial data. However, the geometric-optics construction requires the same
orthogonality relation between $a'(\underline{u})$ and $\nabla_x\phi$ to
hold for positive times. The purpose of this subsubsection is therefore to
show that this condition is preserved by the evolution of the smooth
ground-state and the phase.

\medskip
\noindent
This propagation property is the first key ingredient in the proof of
Theorem~\ref{propagation-thm-intro}: it will allow the corresponding
order-one term in the WKB expansion to vanish throughout the domain
$\Omega$. The same transport computation also yields two related properties
that will be useful later: the propagation of analogous orthogonality
conditions for a general vector field $G$, and the propagation of a uniform
violation of the compatibility condition. These properties are collected in
the following lemma.

\begin{lem}[Propagation of compatibility conditions]
\label{Lemma propagation}
The following properties hold.

\medskip
\noindent
\textbf{(i) Propagation of the compatibility condition.}

\medskip
\noindent
If $a'(\underline u_0)\cdot\nabla_x\phi_0=0$ on $\Omega_0$, then
$a'(\underline u)\cdot\nabla_x\phi=0$ on $\Omega$.

\medskip
\noindent
\textbf{(ii) Propagation for a general $G$.}

\medskip
\noindent
Assume {\rm (i)} and let $G\in C^1(\mathbb R;\mathbb R^d)$. If
$G(\underline u_0)\cdot\nabla_x\phi_0=0$ on $\Omega_0$, then
$G(\underline u)\cdot\nabla_x\phi=0$ on $\Omega$.

\medskip
\noindent
\textbf{(iii) Propagation of uniform violation.}

\medskip
\noindent
If
$\displaystyle
\inf_{x\in\Omega_0}
\left|a'(\underline u_0(x))\cdot\nabla_x\phi_0(x)\right|>0$,
then
$\displaystyle
\inf_{(t,x)\in\Omega}
\left|a'(\underline u(t,x))\cdot\nabla_x\phi(t,x)\right|>0$.

\end{lem}

\begin{Proof}
Set $X:=\partial_t+a(\underline u)\cdot\nabla_x$. Since $\underline u$ and $\phi$ solve their transport equations, $X\underline u=X\phi=0$.

\medskip
\noindent
Let $G\in C^1(\mathbb R;\mathbb R^d)$ and set
$g:=G(\underline u)\cdot\nabla_x\phi$. Differentiating $X\phi=0$ with respect to $x_k$ gives
\[
X(\partial_k\phi)
=
-\partial_k\underline u\,
\bigl(a'(\underline u)\cdot\nabla_x\phi\bigr).
\]
Since $X(G_k(\underline u))=0$, we obtain
\begin{equation}\label{Transport identity}
Xg
=
-\bigl(G(\underline u)\cdot\nabla_x\underline u\bigr)
\bigl(a'(\underline u)\cdot\nabla_x\phi\bigr).
\end{equation}

\medskip
\noindent
For {\rm (i)}, take $G=a'$, hence $g=a'(\underline u)\cdot\nabla_x\phi$. Then \cref{Transport identity} becomes
\begin{equation}\label{Transport divergence}
Xg=-\bigl(\operatorname{div}_x a(\underline u)\bigr)g.
\end{equation}
Since $g=0$ on $\Omega_0$ by \cref{eq:compat-p2}, uniqueness along the characteristics of $X$ gives $g\equiv0$ on $\Omega$.

\medskip
\noindent
For {\rm (ii)}, \cref{Transport identity} together with {\rm (i)} gives $Xg=0$. Since $g=0$ on $\Omega_0$, we again obtain $g\equiv0$ on $\Omega$.

\medskip
\noindent
For {\rm (iii)}, let $g:=a'(\underline u)\cdot\nabla_x\phi$. Along a characteristic $t\mapsto x(t)$ of $X$, \cref{Transport divergence} gives
\[
g(t,x(t))
=
g(0,x(0))
\exp\!\left(
-\int_0^t
\operatorname{div}_x a(\underline u)(s,x(s))\,ds
\right).
\]
Since $\underline u$ is $C^1$ on the compact set $\Omega$, there exists $M>0$ such that
$|\operatorname{div}_x a(\underline u)|\le M$ on $\Omega$. Hence
\[
|g(t,x(t))|
\ge
|g(0,x(0))|e^{-MT},
\]
and therefore
\[
\inf_{(t,x)\in\Omega}|g(t,x)|
\ge
e^{-MT}\inf_{x\in\Omega_0}|g(0,x)|>0.
\]

\medskip
\noindent
And we conclude the proof of this Lemma
\end{Proof}
\subsubsection{Geometric optics reduction}
\label{geometric-optics-reduction}

\medskip
\noindent
The purpose of this subsubsection is to obtain the asymptotic description of
the oscillating solution. We first perform the WKB expansion and identify the
equations satisfied by the phase and the limiting profile. We then introduce
an exact $\varepsilon$-dependent profile and compare it with the limiting one
by the classical stability theory for the associated characteristic systems.

\paragraph{WKB expansion.}\label{WKB expansion}
We begin with a WKB calculus. The main point is to track the successive
powers of $\varepsilon$ and identify the terms that vanish thanks to the
equations satisfied by the ground-state and the phase, together with the
propagated compatibility condition. Given a general profile
$U(t,x,\theta)$, we seek a solution in the form
\begin{equation}
v^{\varepsilon}(t,x)
:=
\underline{u}(t,x)
+
\varepsilon
U\left(
t,x,\frac{\phi(t,x)}{\varepsilon^{2}}
\right).
\end{equation}

\medskip
\noindent
To make the remainder explicit, we use the Taylor expansion
\begin{equation}\label{eq:Taylor-F-remainder}
\begin{aligned}
F(\underline u+\varepsilon U)
={}&
F(\underline u)
+\varepsilon a(\underline u)U
+\frac{\varepsilon^2}{2}a'(\underline u)U^2
+\frac{\varepsilon^3}{6}a''(\underline u)U^3
+\varepsilon^4\widetilde R^\varepsilon,
\end{aligned}
\end{equation}
where
\begin{equation}\label{eq:def-Rtilde-eps}
\widetilde R^\varepsilon(t,x,\theta)
:=
\frac{U(t,x,\theta)^4}{6}
\int_0^1
(1-s)^3
a'''\!\left(
\underline u(t,x)
+s\varepsilon U(t,x,\theta)
\right)\,ds .
\end{equation}

\medskip
\noindent
We define the $(d+1)$-dimensional profile remainder by
\begin{equation}\label{eq:def-R-eps}
R^\varepsilon(t,x,\theta)
:=
\left(
\frac12 a'(\underline u)U^2
+
\frac{\varepsilon}{6}a''(\underline u)U^3
+
\varepsilon^2\widetilde R^\varepsilon
\,;\,
\nabla_x\phi\cdot\widetilde R^\varepsilon
\right),
\end{equation}
and the physical remainder $r^\varepsilon$ through
\begin{equation}\label{eq:def-r-eps}
\operatorname{div}_x r^\varepsilon(t,x)
:=
\left.
\operatorname{div}_{(x,\theta)}
R^\varepsilon(t,x,\theta)
\right|_{\theta=\phi(t,x)/\varepsilon^2}.
\end{equation}

\medskip
\noindent
A formal WKB computation as in \cite{J2,J4} gives
\begin{equation}\label{WKB-expansion}
\begin{aligned}
\partial_t v^{\varepsilon}+\operatorname{div}_x F(v^{\varepsilon})
={}&\varepsilon^{-1}
\Big(\partial_t\phi+a(\underline{u})\cdot\nabla_x\phi\Big)
\,\partial_\theta U \\
&+\Big(\partial_t \underline{u}
+\operatorname{div}_x F(\underline{u})\Big)
+\frac12\Big(a'(\underline{u})\cdot\nabla_x\phi\Big)
\,\partial_\theta \big(U^{2}\big) \\
&+\varepsilon\Big(
\partial_t U+\operatorname{div}_x\!\big(U\,a(\underline{u})\big)
+\frac16\big(a''(\underline{u})\cdot\nabla_x\phi\big)
\,\partial_\theta \big(U^{3}\big)
\Big) \\
&+\varepsilon^{2}\,
\operatorname{div}_x\big(r^{\varepsilon}(t,x)\big).
\end{aligned}
\end{equation}
Here the profiles $\partial_\theta(U^j)$ are evaluated at
$\theta=\phi(t,x)/\varepsilon^2$. The remainders
$\widetilde R^\varepsilon$ and $R^\varepsilon$ are smooth and uniformly
bounded, together with the derivatives involved above, for
$\varepsilon\in(0,1]$ on the considered compact domain.

\medskip
\noindent
In \cref{WKB-expansion}, the term of order $\varepsilon^{-1}$ vanishes by
the eikonal equation \cref{eq:eikonal}. At order $\varepsilon^0$, the
ground-state term vanishes by \cref{eq:claw}, while the quadratic term
vanishes by the compatibility condition, propagated throughout $\Omega$ by
Lemma~\ref{Lemma propagation}. Finally, the term of order $\varepsilon$
vanishes when $U$ solves the profile equation \cref{Profile-intro}. Hence
\begin{equation}\label{approx-conservation-v}
\partial_t v^\varepsilon+\operatorname{div}_x F(v^\varepsilon)
=
\varepsilon^{2}\,\operatorname{div}_x
\big(r^\varepsilon(t,x)\big).
\end{equation}

\paragraph{Exact $\varepsilon$-dependent profile.}
\label{exact-profile-reduction}
Rather than pursuing the WKB expansion to higher orders, we follow the
strategy of \cite{J2,J4} and introduce an exact $\varepsilon$-dependent
profile $U^\varepsilon$. The corresponding reconstruction satisfies the
original conservation law exactly, reducing the problem to a comparison
between $U^\varepsilon$ and the limiting profile $U$.

\medskip
\noindent
More precisely, $U^\varepsilon$ is introduced through
\begin{equation}\label{def-ueps}
u^\varepsilon(t,x):=\underline u(t,x)
+\varepsilon\,U^\varepsilon\!\left(
t,x,\frac{\phi(t,x)}{\varepsilon^{2}}
\right),
\end{equation}
with $U^\varepsilon(0,x,\theta)=U_0(x,\theta)$, so that
\cref{eq:claw,eq:initdata} is satisfied.

\medskip
\noindent
Set $y=(x,\theta)\in\Omega_0\times\mathbb T$ and let
$\operatorname{div}_y$ denote the divergence in $(x,\theta)$. The same WKB
computation yields
\begin{equation}\label{CL-Ueps}
\partial_t U^\varepsilon
+
\operatorname{div}_y
G^\varepsilon\!\bigl(t,x,U^\varepsilon\bigr)
=
0,
\qquad
U^\varepsilon(0,x,\theta)=U_0(x,\theta),
\end{equation}
while the limiting profile $U$ satisfies
\begin{equation}\label{CL-U}
\partial_t U
+
\operatorname{div}_y
G\!\bigl(t,x,U\bigr)
=
0,
\qquad
U(0,x,\theta)=U_0(x,\theta).
\end{equation}
The corresponding fluxes satisfy
\begin{equation}\label{Geps-G}
G^\varepsilon(t,x,V)
=
G(t,x,V)
+
\varepsilon R^\varepsilon(t,x,V),
\end{equation}
where
\begin{equation}\label{def-G}
G(t,x,V)
=
\Bigl(
V\,a(\underline u(t,x))
\; ;\;
\frac16
\bigl(
a''(\underline u(t,x))
\cdot
\nabla_x\phi(t,x)
\bigr)V^3
\Bigr).
\end{equation}
Here $R^\varepsilon(t,x,V)$ is given by the same expression as in
\eqref{eq:def-R-eps}, with $U(t,x,\theta)$ replaced by the generic variable
$V$ and $\widetilde R^\varepsilon$ correspondingly given by
\eqref{eq:def-Rtilde-eps}. Thus $R^\varepsilon$ is smooth and uniformly
bounded on $\Omega\times\{|V|\leq C\}$, uniformly for
$\varepsilon\in(0,1]$, for every fixed $C>0$.

\medskip
\noindent
Since the two profile equations have the same smooth initial datum and
$G^\varepsilon=G+\varepsilon R^\varepsilon$, the classical local theory and
smooth dependence on the parameter $\varepsilon$ for the associated
characteristic systems \cite{PontryaginODE} yield, after reducing $T>0$ if necessary, a common
smooth existence interval and
\begin{equation}\label{eq:profile-stability}
\sup_{t\in[0,T]}
\left\|
U^\varepsilon(t,\cdot)-U(t,\cdot)
\right\|_{L^\infty(\Omega_0\times\mathbb T)}
\leq
C_T\varepsilon ,
\end{equation}
uniformly for $0<\varepsilon\leq1$.

\medskip
\noindent
Using \cref{def-ueps} and the definition of $v^\varepsilon$, we therefore
obtain
\[
u^\varepsilon(t,x)-v^\varepsilon(t,x)
=
\varepsilon
\left(
U^\varepsilon-U
\right)
\left(
t,x,\frac{\phi(t,x)}{\varepsilon^2}
\right)
=
O(\varepsilon^2)
\]
uniformly on $\Omega$. Hence
\begin{equation}\label{eq:propagation-final-error}
\left\|
u^\varepsilon(t,\cdot)
-
\underline u(t,\cdot)
-
\varepsilon
U\left(
t,\cdot,\frac{\phi(t,\cdot)}{\varepsilon^2}
\right)
\right\|_{C^0(\Omega_t)}
\leq
C_T\varepsilon^2,
\qquad t\in[0,T].
\end{equation}

\medskip
\noindent
Finally, since $\underline u$, $\phi$, and $U^\varepsilon$ are smooth on the
common interval $[0,T]$, the reconstruction \cref{def-ueps} defines a smooth
solution of \cref{eq:claw,eq:initdata} on $\Omega$. By uniqueness of
classical solutions, it coincides with the maximal classical solution issued
from $u_0^\varepsilon$ on their common interval of existence. Therefore, if
$T_\varepsilon$ denotes its maximal lifespan,
\[
T_\varepsilon\geq T,
\qquad
0<\varepsilon\leq1.
\]
Together with \cref{eq:propagation-final-error}, this proves
Theorem~\ref{propagation-thm-intro}.

\subsection{Cancellation}\label{Cancellation}
We now turn to the cancellation of highly oscillating solutions in the case
$p=2$. As stated in Theorem~\ref{cancellation-thm-intro}, when the
compatibility condition \cref{eq:compat-p2} is uniformly violated, in the
sense of \cref{violation-intro}, the oscillations are cancelled and only the
mean of the initial oscillating profile is propagated, up to a remainder.
In contrast with the propagation regime, shock formation occurs on a short
time scale, so that weak entropy solutions in the sense of
Kružkov~\cite{Kr} have to be considered, as in \cite{CJR,J4}.

\medskip
\noindent
For $k\in\mathbb Z^d$, we introduce the rescaled functions
\begin{equation}\label{eq:local-rescaling}
    v_k^\varepsilon(\tau,y)
    :=
    u^\varepsilon(\varepsilon\tau,\varepsilon k+\varepsilon y),
    \qquad
    \underline v_k^\varepsilon(\tau,y)
    :=
    \underline u(\varepsilon\tau,\varepsilon k+\varepsilon y),
\end{equation}
and
\begin{equation}\label{eq:local-averaged-profile}
    \overline U_k^\varepsilon(\tau,y)
    :=
    \overline U(\varepsilon\tau,\varepsilon k+\varepsilon y).
\end{equation}

\medskip
\noindent The initial decay of the oscillations is analyzed in the time boundary layer
$0\leq t\leq\varepsilon\mathcal I$, where $\mathcal I$ is chosen
sufficiently large but independent of $\varepsilon$. At the end of this
boundary layer, the oscillating contribution is sufficiently small for the
resulting estimate to be propagated to later times by the $L^1$-contraction
principle. The corresponding localization is illustrated in
Figure~\ref{fig:boundary-layer}.

\begin{figure}[H]
    \centering
    \includegraphics[width=0.55\textwidth]{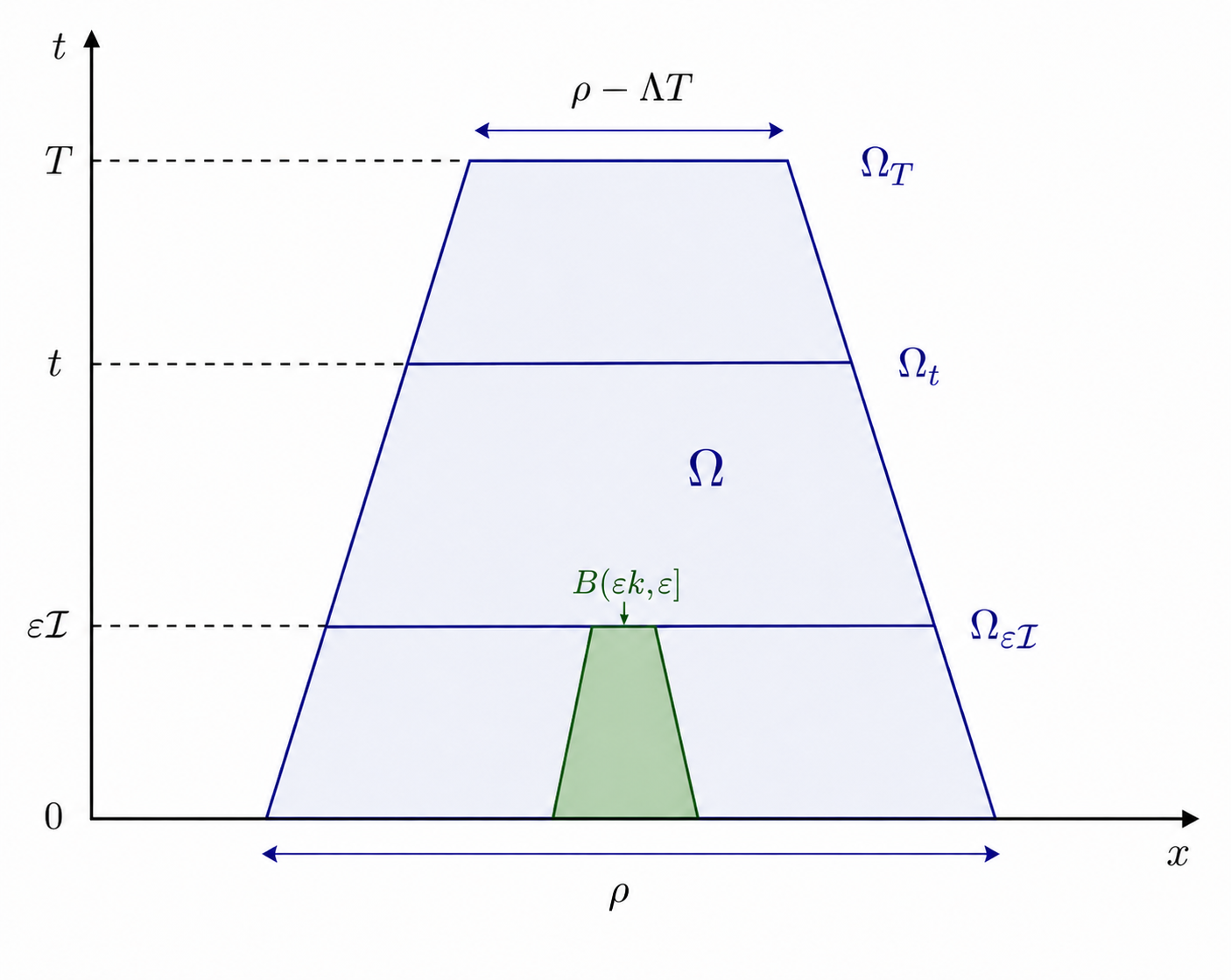}
    \caption{Spacetime domain $\Omega$ and localization in the initial time
    layer $0\le t\le \varepsilon\mathcal I$.}
    \label{fig:boundary-layer}
\end{figure}

\medskip
\noindent The key estimate obtained at the end of this boundary layer is the
following.

\begin{prop}[Local cancellation estimate]\label{prop:local-cancellation}
Let $\rho>0$. There exist a time $\mathcal I_0>1$ and a constant $C_\rho>0$
such that, for every $\mathcal I\geq\mathcal I_0$ with
$\varepsilon\mathcal I\leq T$, there exists
$\varepsilon_{\mathcal I}>0$ such that, for every
$0<\varepsilon\leq\varepsilon_{\mathcal I}$ and every
$k\in\mathbb Z^d$ satisfying
\[
    \varepsilon k\in B(0,\rho],
\]
one has
\begin{equation}\label{eq:local-cancellation-estimate}
    \left\|
    v_k^\varepsilon(\mathcal I,\cdot)
    -
    \underline v_k^\varepsilon(\mathcal I,\cdot)
    -
    \varepsilon\overline U_k^\varepsilon(\mathcal I,\cdot)
    \right\|_{L^1(B(0,1])}
    \leq
    C_\rho\frac{\varepsilon}{\mathcal I}.
\end{equation}
\end{prop}

\medskip
\noindent The proof is organized as follows. We first assume
Proposition~\ref{prop:local-cancellation}. In
Subsubsection~\ref{decompoballs}, the localized estimates
\cref{eq:local-cancellation-estimate}, indexed by $k$, are patched together
and transferred from the rescaled variables to the physical variables.
In Subsubsection~\ref{contractionprinciple}, the resulting estimate is then
propagated from the intermediate time $t=\varepsilon\mathcal I$ by means of
the local $L^1$-contraction principle, which allows us to conclude the proof
of Theorem~\ref{cancellation-thm-intro}. We then turn to the proof of
Proposition~\ref{prop:local-cancellation}. In
Subsubsection~\ref{boundarylayer}, the solution is analyzed inside the
initial time boundary layer by means of weakly nonlinear geometric optics (WNGO),
leading to an Oleinik decay estimate for the localized oscillating profile.
Finally, in Subsubsection~\ref{compaave}, this profile is compared with its
average, which completes the proof of
Proposition~\ref{prop:local-cancellation}.

\subsubsection{Decomposition of $\Omega_{\varepsilon\mathcal I}$ into
$\varepsilon$-balls and reconstruction}\label{decompoballs}
We begin by rewriting the local estimate
\cref{eq:local-cancellation-estimate} in the original spatial variables.
We set $x=\varepsilon k+\varepsilon y$. Here and in the sequel, we use the notation
$B(\varepsilon k,\varepsilon]:=\varepsilon k+\varepsilon B(0,1]$.
Thus $B(\varepsilon k,\varepsilon]$ is a hypercube of side length
$\varepsilon$. Moreover,
$y\in B(0,1]$ if and only if $x\in B(\varepsilon k,\varepsilon]$, and
$dx=\varepsilon^d\,dy$.

\medskip
\noindent
Hence \eqref{eq:local-cancellation-estimate} becomes
\[
\frac{1}{\varepsilon^d}
\int_{B(\varepsilon k,\varepsilon]}
\left|
u^\varepsilon(\varepsilon\mathcal I,x)
-
\underline u(\varepsilon\mathcal I,x)
-
\varepsilon \overline U(\varepsilon\mathcal I,x)
\right|\,dx
\leq
C_\rho\frac{\varepsilon}{\mathcal I}.
\]
Multiplying by $\varepsilon^d$, we obtain the estimate on one physical
hypercube:
\begin{equation}\label{eq:one-cube-estimate}
    \left\|
    u^\varepsilon(\varepsilon\mathcal I,\cdot)
    -
    \underline u(\varepsilon\mathcal I,\cdot)
    -
    \varepsilon \overline U(\varepsilon\mathcal I,\cdot)
    \right\|_{L^1(B(\varepsilon k,\varepsilon])}
    \leq
    C_\rho\frac{\varepsilon^{d+1}}{\mathcal I}.
\end{equation}

\medskip
\noindent
We sum \eqref{eq:one-cube-estimate} over the small hypercubes covering
$B(0,\rho]$. By finite speed of propagation, the values of a solution on
$\Omega_t$ at time $t$ can only depend on its values at the earlier time
$\varepsilon\mathcal I$ in the backward domain
$\Omega_t+B(0,\Lambda(t-\varepsilon\mathcal I)]$.
Since $\Omega_t=B(0,\rho-\Lambda t]$, we have
\[
\Omega_t+B(0,\Lambda(t-\varepsilon\mathcal I)]
\subset
B(0,\rho-\Lambda\varepsilon\mathcal I]
\subset
B(0,\rho].
\]
Thus the whole backward domain of dependence of $\Omega_t$ between the
times $\varepsilon\mathcal I$ and $t$ is contained in $B(0,\rho]$.

\medskip
\noindent
We cover $B(0,\rho]$ by the hypercubes
$B(\varepsilon k,\varepsilon]$ with $\varepsilon k\in B(0,\rho]$. Thus,
\begin{equation}\label{eq:covering-BM}
    B(0,\rho]
    \subset
    \bigcup_{\varepsilon k\in B(0,\rho]}
    B(\varepsilon k,\varepsilon].
\end{equation}
Moreover, the number of such hypercubes satisfies
\begin{equation}\label{eq:number-cubes}
    \#\{k\in\mathbb Z^d:\ \varepsilon k\in B(0,\rho]\}
    \leq
    C_\rho\varepsilon^{-d},
\end{equation}
where $C_\rho>0$ is independent of $\varepsilon$.

\medskip
\noindent
Using \eqref{eq:one-cube-estimate} and \eqref{eq:covering-BM}, we obtain
\[
\begin{aligned}
&\left\|
u^\varepsilon(\varepsilon\mathcal I,\cdot)
-
\underline u(\varepsilon\mathcal I,\cdot)
-
\varepsilon \overline U(\varepsilon\mathcal I,\cdot)
\right\|_{L^1(B(0,\rho])}
\\
&\leq
\sum_{\varepsilon k\in B(0,\rho]}
\left\|
u^\varepsilon(\varepsilon\mathcal I,\cdot)
-
\underline u(\varepsilon\mathcal I,\cdot)
-
\varepsilon \overline U(\varepsilon\mathcal I,\cdot)
\right\|_{L^1(B(\varepsilon k,\varepsilon])}
\\
&\leq
\sum_{\varepsilon k\in B(0,\rho]}
C_\rho\frac{\varepsilon^{d+1}}{\mathcal I}.
\end{aligned}
\]
By \eqref{eq:number-cubes},
\[
\sum_{\varepsilon k\in B(0,\rho]}
C_\rho\frac{\varepsilon^{d+1}}{\mathcal I}
\leq
C_\rho\varepsilon^{-d}
C_\rho\frac{\varepsilon^{d+1}}{\mathcal I}
\leq
C_\rho\frac{\varepsilon}{\mathcal I},
\]
after enlarging $C_\rho$ if necessary. Therefore,
\begin{equation}\label{eq:estimate-intermediate-time}
    \left\|
    u^\varepsilon(\varepsilon\mathcal I,\cdot)
    -
    \underline u(\varepsilon\mathcal I,\cdot)
    -
    \varepsilon \overline U(\varepsilon\mathcal I,\cdot)
    \right\|_{L^1(B(0,\rho])}
    \leq
    C_\rho\frac{\varepsilon}{\mathcal I}.
\end{equation}

\subsubsection{$L^1$-contraction principle}
\label{contractionprinciple}

We now propagate the estimate \eqref{eq:estimate-intermediate-time} from the
intermediate time $t=\varepsilon\mathcal I$ to the final time $t$.
Since $\underline u+\varepsilon\overline U$ is not an exact solution of the
nonlinear conservation law, we cannot directly apply the $L^1$-contraction
principle between $u^\varepsilon$ and
$\underline u+\varepsilon\overline U$. We therefore introduce an exact
comparison solution.

\medskip
\noindent
Let $\widehat u^\varepsilon$ be the entropy solution of
\begin{equation}\label{eq:uhat-equation}
    \partial_s \widehat u^\varepsilon
    +
    \operatorname{div}_x F(\widehat u^\varepsilon)
    =
    0,
    \qquad
    \varepsilon\mathcal I<s\leq t,
\end{equation}
with initial datum at time $s=\varepsilon\mathcal I$ given by
\begin{equation}\label{eq:uhat-initial-data}
    \widehat u^\varepsilon(\varepsilon\mathcal I,x)
    =
    \underline u(\varepsilon\mathcal I,x)
    +
    \varepsilon\overline U(\varepsilon\mathcal I,x).
\end{equation}
Then, by \eqref{eq:estimate-intermediate-time} and \eqref{eq:uhat-initial-data},
\begin{equation}\label{eq:u-uhat-intermediate}
    \left\|
    u^\varepsilon(\varepsilon\mathcal I,\cdot)
    -
    \widehat u^\varepsilon(\varepsilon\mathcal I,\cdot)
    \right\|_{L^1(B(0,\rho])}
    \leq
    C_\rho\frac{\varepsilon}{\mathcal I}.
\end{equation}

\medskip
\noindent
Now $u^\varepsilon$ and $\widehat u^\varepsilon$ are entropy solutions of the
same scalar conservation law on the time interval
$[\varepsilon\mathcal I,t]$. Hence, by the local $L^1$-contraction property,
\begin{equation}\label{eq:L1-contraction-backward-cone}
\begin{aligned}
&\left\|
u^\varepsilon(t,\cdot)
-
\widehat u^\varepsilon(t,\cdot)
\right\|_{L^1(\Omega_t)}
\\
&\leq
\left\|
u^\varepsilon(\varepsilon\mathcal I,\cdot)
-
\widehat u^\varepsilon(\varepsilon\mathcal I,\cdot)
\right\|_{L^1\left(\Omega_t+B(0,\Lambda(t-\varepsilon\mathcal I)]\right)}.
\end{aligned}
\end{equation}
The set on the right-hand side is contained in $B(0,\rho]$. Therefore,
\begin{equation}\label{eq:L1-contraction-final}
\begin{aligned}
&\left\|
u^\varepsilon(t,\cdot)
-
\widehat u^\varepsilon(t,\cdot)
\right\|_{L^1(\Omega_t)}
\\
&\leq
\left\|
u^\varepsilon(\varepsilon\mathcal I,\cdot)
-
\widehat u^\varepsilon(\varepsilon\mathcal I,\cdot)
\right\|_{L^1(B(0,\rho])}.
\end{aligned}
\end{equation}
Combining \eqref{eq:u-uhat-intermediate} and \eqref{eq:L1-contraction-final},
we obtain
\begin{equation}\label{eq:u-uhat-final-estimate}
    \left\|
    u^\varepsilon(t,\cdot)
    -
    \widehat u^\varepsilon(t,\cdot)
    \right\|_{L^1(\Omega_t)}
    \leq
    C_\rho\frac{\varepsilon}{\mathcal I}.
\end{equation}

\medskip
\noindent
It remains to compare the exact comparison solution $\widehat u^\varepsilon$
with the first order expansion $\underline u+\varepsilon\overline U$.
This comparison is classical, because we are in the smooth regime. Indeed,
$\underline u$ is smooth on $[0,T]$ by assumption, and $\overline U$ solves
the linearized equation \eqref{Profile-2-intro} with smooth coefficients.
Hence $\overline U$ is also smooth on $[0,T]$, and therefore
$\underline u+\varepsilon\overline U$ is smooth.

\medskip
\noindent
Set $w^\varepsilon:=\underline u+\varepsilon\overline U$. By the definition of
$\widehat u^\varepsilon$, we have
\begin{equation}\label{eq:uhat-w-same-data}
    \widehat u^\varepsilon(\varepsilon\mathcal I,\cdot)
    =
    w^\varepsilon(\varepsilon\mathcal I,\cdot).
\end{equation}
Moreover, $w^\varepsilon$ is an approximate solution of the nonlinear
conservation law up to an error of order $\varepsilon^2$. Indeed,
\[
\begin{aligned}
\partial_s w^\varepsilon+\operatorname{div}_x F(w^\varepsilon)
&=
\partial_s(\underline u+\varepsilon\overline U)
+
\operatorname{div}_x F(\underline u+\varepsilon\overline U).
\end{aligned}
\]
Using Taylor expansion of $F$ around $\underline u$, we get
\[
F(\underline u+\varepsilon\overline U)
=
F(\underline u)
+
\varepsilon F'(\underline u)\overline U
+
O(\varepsilon^2).
\]
Hence
\[
\begin{aligned}
\partial_s w^\varepsilon+\operatorname{div}_x F(w^\varepsilon)
&=
\underbrace{
\partial_s\underline u
+
\operatorname{div}_x F(\underline u)
}_{=0}
\\
&\quad+
\varepsilon
\underbrace{
\left(
\partial_s\overline U
+
\operatorname{div}_x\big(F'(\underline u)\overline U\big)
\right)
}_{=0}
+
O(\varepsilon^2).
\end{aligned}
\]
Therefore
\begin{equation}\label{eq:w-approx-solution}
    \partial_s w^\varepsilon
    +
    \operatorname{div}_x F(w^\varepsilon)
    =
    O(\varepsilon^2).
\end{equation}

\medskip
\noindent
Since $\widehat u^\varepsilon$ is the exact solution of
$\partial_s\widehat u^\varepsilon+\operatorname{div}_xF(\widehat u^\varepsilon)=0$
with the same initial datum as $w^\varepsilon$ at time
$\varepsilon\mathcal I$, and since $w^\varepsilon$ satisfies the same
equation up to the residual \eqref{eq:w-approx-solution}, the classical stability
estimate for smooth solutions gives
\begin{equation}\label{eq:uhat-w-comparison}
    \left\|
    \widehat u^\varepsilon(t,\cdot)
    -
    w^\varepsilon(t,\cdot)
    \right\|_{L^1(\Omega_t)}
    \leq
    C_t\varepsilon^2.
\end{equation}
Equivalently,
\begin{equation}\label{eq:uhat-background-comparison}
    \left\|
    \widehat u^\varepsilon(t,\cdot)
    -
    \underline u(t,\cdot)
    -
    \varepsilon\overline U(t,\cdot)
    \right\|_{L^1(\Omega_t)}
    \leq
    C_t\varepsilon^2.
\end{equation}

\medskip
\noindent
We can now conclude by combining the two estimates. By the triangle inequality,
\[
\begin{aligned}
&\left\|
u^\varepsilon(t,\cdot)
-
\underline u(t,\cdot)
-
\varepsilon\overline U(t,\cdot)
\right\|_{L^1(\Omega_t)}
\\
&\leq
\left\|
u^\varepsilon(t,\cdot)
-
\widehat u^\varepsilon(t,\cdot)
\right\|_{L^1(\Omega_t)}
+
\left\|
\widehat u^\varepsilon(t,\cdot)
-
\underline u(t,\cdot)
-
\varepsilon\overline U(t,\cdot)
\right\|_{L^1(\Omega_t)}.
\end{aligned}
\]
Using \eqref{eq:u-uhat-final-estimate} and
\eqref{eq:uhat-background-comparison}, we obtain
\begin{equation}\label{eq:final-before-choice}
    \left\|
    u^\varepsilon(t,\cdot)
    -
    \underline u(t,\cdot)
    -
    \varepsilon\overline U(t,\cdot)
    \right\|_{L^1(\Omega_t)}
    \leq
    C_\rho\frac{\varepsilon}{\mathcal I}
    +
    C_t\varepsilon^2.
\end{equation}

\medskip
\noindent
It remains to choose the parameters. Let $\ell>0$ be fixed. Since
Proposition~\ref{prop:local-cancellation} holds for every sufficiently large
fast time $\mathcal I$, we choose $\mathcal I$ large enough so that
$\mathcal I\geq\mathcal I_0$ and $C_\rho/\mathcal I\leq\ell/2$.
With this $\mathcal I$ fixed, we choose $\varepsilon_\ell>0$ small enough so
that, for every $0<\varepsilon\leq\varepsilon_\ell$, the local proposition
applies, namely $0<\varepsilon\leq\varepsilon_{\mathcal I}$, and also
$\varepsilon\mathcal I<t$ and $C_t\varepsilon\leq\ell/2$.
Therefore, from \eqref{eq:final-before-choice}, we get
\[
\begin{aligned}
\left\|
u^\varepsilon(t,\cdot)
-
\underline u(t,\cdot)
-
\varepsilon\overline U(t,\cdot)
\right\|_{L^1(\Omega_t)}
&\leq
\frac{\ell}{2}\varepsilon
+
\frac{\ell}{2}\varepsilon
=
\ell\varepsilon.
\end{aligned}
\]
Since $\ell>0$ is arbitrary, we conclude that
\[
\left\|
u^\varepsilon(t,\cdot)
-
\underline u(t,\cdot)
-
\varepsilon\overline U(t,\cdot)
\right\|_{L^1(\Omega_t)}
=
o(\varepsilon).
\]

\medskip
\noindent
There remains now to prove Proposition~\ref{prop:local-cancellation}.

\subsubsection{Analysis in the boundary layer}\label{boundarylayer}
Some substantial preparatory work needs to be carried out near $t=0$, in a
time boundary layer of size $\varepsilon$. This is achieved in this subsection
which contains core arguments. In Paragraph \ref{weaklynonlinear}, the section
$\Omega_{\varepsilon\mathcal I}$ is partitioned into small balls of radius
$\varepsilon$. By this way, the discussion can fall within the framework of
weakly nonlinear geometric optics. The rescaled solution can then be compared
to some WKB ansatz, with a profile which satisfies a law of Burger's type.
As explained in Paragraph \ref{paraOleinik}, it becomes possible to establish some Oleinik
estimate, leading to the supnorm decay of the profile. In Paragraph \ref{paraOleinik},
this is the indispensable key to be able to compare the WKB expansion and the
averaged model.

\paragraph{Weakly nonlinear geometric optics}\label{weaklynonlinear} 

\medskip
\noindent
We localize the problem on hypercubes of size $\varepsilon$. For
$k\in\mathbb Z^d$, we introduce the rescaled variables
$t=\varepsilon\tau$, $x=\varepsilon k+\varepsilon y$, and define
\[
v_k^\varepsilon(\tau,y)
:=
u^\varepsilon(\varepsilon\tau,\varepsilon k+\varepsilon y),
\qquad
\underline v_k^\varepsilon(\tau,y)
:=
\underline u(\varepsilon\tau,\varepsilon k+\varepsilon y).
\]

\medskip
\noindent
Since $u^\varepsilon$ is an entropy solution of \cref{eq:claw},
\cref{eq:initdata}, for every constant $\kappa\in\mathbb R$ it satisfies
the Kruzkov entropy inequality
\[
\partial_t |u^\varepsilon-\kappa|
+
\operatorname{div}_x
\left[
\operatorname{sgn}(u^\varepsilon-\kappa)
\bigl(F(u^\varepsilon)-F(\kappa)\bigr)
\right]
\leq 0
\]
in the sense of distributions. To transfer this property to the rescaled
variables, let
$\psi\in C_c^\infty((0,+\infty)\times\mathbb R^d)$, $\psi\geq0$, and set
$\phi(t,x):=\psi(t/\varepsilon,(x-\varepsilon k)/\varepsilon)$.
Then $\partial_t\phi=\varepsilon^{-1}\partial_\tau\psi$,
$\nabla_x\phi=\varepsilon^{-1}\nabla_y\psi$, while
$dt\,dx=\varepsilon^{d+1}\,d\tau\,dy$. Using $\phi$ in the entropy
inequality and changing variables gives
\[
\begin{aligned}
0\leq
\varepsilon^d
\int_0^{+\infty}\int_{\mathbb R^d}
&|v_k^\varepsilon-\kappa|\partial_\tau\psi
\\
&+
\operatorname{sgn}(v_k^\varepsilon-\kappa)
\bigl(F(v_k^\varepsilon)-F(\kappa)\bigr)\cdot\nabla_y\psi
\,dy\,d\tau .
\end{aligned}
\]
Since $\varepsilon^d>0$, $v_k^\varepsilon$ satisfies the corresponding
Kruzkov entropy inequality in $(\tau,y)$, and is therefore an entropy
solution of
\begin{equation}\label{eq:rescaled-local-scl}
\partial_\tau v_k^\varepsilon
+
\operatorname{div}_y F(v_k^\varepsilon)
=
0.
\end{equation}

\medskip
\noindent
We now rewrite the initial data in the rescaled variables. By definition,
$v_k^\varepsilon(0,y)
=u_0^\varepsilon(\varepsilon k+\varepsilon y)$, and \cref{eq:initdata}
gives
\begin{equation}\label{eq:local-initial-data-exact}
v_k^\varepsilon(0,y)
=
\underline u_0(\varepsilon k+\varepsilon y)
+
\varepsilon
U_0\left(
\varepsilon k+\varepsilon y,
\frac{\varphi_0(\varepsilon k+\varepsilon y)}{\varepsilon^2}
\right).
\end{equation}
Since
$\underline v_k^\varepsilon(0,y)
=\underline u_0(\varepsilon k+\varepsilon y)$,
this becomes
\begin{equation}\label{eq:local-initial-data-background}
v_k^\varepsilon(0,y)
=
\underline v_k^\varepsilon(0,y)
+
\varepsilon
U_0\left(
\varepsilon k+\varepsilon y,
\frac{\varphi_0(\varepsilon k+\varepsilon y)}{\varepsilon^2}
\right).
\end{equation}

\medskip
\noindent
We now expand the phase around the point $\varepsilon k$. For
$y\in B(0,1]$, Taylor's formula gives
\[
\varphi_0(\varepsilon k+\varepsilon y)
=
\varphi_0(\varepsilon k)
+
\varepsilon \nabla_x\varphi_0(\varepsilon k)\cdot y
+
\frac{\varepsilon^2}{2}
D_x^2\varphi_0(\varepsilon k)(y,y)
+
\varepsilon^3 R_k^\varepsilon(y),
\]
where $R_k^\varepsilon$ is smooth and uniformly bounded, for
$\varepsilon k$ in a fixed compact set and $y\in B(0,1]$.

\medskip
\noindent
Dividing by
$\varepsilon^2$, we obtain
\begin{equation}\label{eq:local-phase-expansion}
\frac{\varphi_0(\varepsilon k+\varepsilon y)}{\varepsilon^2}
=
\frac{\varphi_0(\varepsilon k)}{\varepsilon^2}
+
\frac{\nabla_x\varphi_0(\varepsilon k)\cdot y}{\varepsilon}
+
\frac12 D_x^2\varphi_0(\varepsilon k)(y,y)
+
\varepsilon R_k^\varepsilon(y).
\end{equation}

\medskip
\noindent
We introduce the local linear phase
\begin{equation}\label{eq:local-linear-phase}
\varphi_k^\varepsilon(y)
:=
\nabla_x\varphi_0(\varepsilon k)\cdot y.
\end{equation}
Then \cref{eq:local-phase-expansion} reads
\[
\frac{\varphi_0(\varepsilon k+\varepsilon y)}{\varepsilon^2}
=
\frac{\varphi_0(\varepsilon k)}{\varepsilon^2}
+
\frac{\varphi_k^\varepsilon(y)}{\varepsilon}
+
\frac12 D_x^2\varphi_0(\varepsilon k)(y,y)
+
\varepsilon R_k^\varepsilon(y).
\]

\medskip
\noindent
Therefore the initial data \cref{eq:local-initial-data-background} can be
written in the exact weakly nonlinear geometric optics form
\begin{equation}\label{eq:local-initial-data-WNGO}
v_k^\varepsilon(0,y)
=
\underline v_k^\varepsilon(0,y)
+
\varepsilon
\mathcal V_{k,0}^\varepsilon
\left(
y,
\frac{\varphi_k^\varepsilon(y)}{\varepsilon}
\right),
\end{equation}
where
\begin{equation}\label{eq:local-initial-profile}
\mathcal V_{k,0}^\varepsilon(y,\theta)
:=
U_0\left(
\varepsilon k+\varepsilon y,
\frac{\varphi_0(\varepsilon k)}{\varepsilon^2}
+
\theta
+
\frac12 D_x^2\varphi_0(\varepsilon k)(y,y)
+
\varepsilon R_k^\varepsilon(y)
\right).
\end{equation}
\medskip
\noindent
Hence the localized problem is exactly in the standard weakly nonlinear
geometric optics regime.

\medskip
\noindent
We can now apply the nonlinear geometric optics theorem as in~\cite{Ch92,Ch97}.
Indeed, \cref{eq:rescaled-local-scl} has initial data of the standard weakly
nonlinear geometric optics form \cref{eq:local-initial-data-WNGO}, where
$\mathcal V_{k,0}^\varepsilon$ is smooth and $1$-periodic with respect to the
second variable, and where the local phase $\varphi_k^\varepsilon$ defined in
\cref{eq:local-linear-phase} is smooth and linear.

\medskip
\noindent
The dependence of $\mathcal V_{k,0}^\varepsilon$ and
$\varphi_k^\varepsilon$ on $\varepsilon$ and $k$ is harmless. Indeed,
$\varepsilon k$ remains in a fixed compact set, the phase shift
$\varphi_0(\varepsilon k)/\varepsilon^2$ is only a translation in the
periodic variable, and all the profiles are uniformly smooth and uniformly
bounded with respect to $(k,\varepsilon)$.

\medskip
\noindent
Therefore, by~\cite{Ch92,Ch97}, for every fixed $\mathcal I>0$, we have
\begin{equation}\label{eq:local-WNGO-expansion}
\left\|
v_k^\varepsilon(\tau,\cdot)
-
\underline v_k^\varepsilon(\tau,\cdot)
-
\varepsilon
\mathcal V_k^\varepsilon
\left(
\tau,\cdot,
\frac{\widehat\varphi_k^\varepsilon(\tau,\cdot)}{\varepsilon}
\right)
\right\|_{L^1(B(0,1])}
\leq
C_{\mathcal I}\varepsilon^2,
\end{equation}
for $0\leq \tau\leq \mathcal I$, where the constant $C_{\mathcal I}$ is
independent of $k$ and $\varepsilon$.

\medskip
\noindent
Here $\widehat\varphi_k^\varepsilon$ is the solution of the eikonal equation
\begin{equation}\label{eq:local-eikonal}
\begin{cases}
\partial_\tau \widehat\varphi_k^\varepsilon
+
a(\underline v_k^\varepsilon)
\cdot
\nabla_y\widehat\varphi_k^\varepsilon
=
0,
\\[1mm]
\widehat\varphi_k^\varepsilon(0,y)
=
\varphi_k^\varepsilon(y).
\end{cases}
\end{equation}
Moreover, $\mathcal V_k^\varepsilon$ is the corresponding oscillating
profile. It solves the Burgers-type profile equation
\begin{equation}\label{eq:local-profile-equation}
\begin{cases}
\partial_\tau \mathcal V_k^\varepsilon
+
\operatorname{div}_y
\left(
a(\underline v_k^\varepsilon)\mathcal V_k^\varepsilon
\right)
+
\partial_\theta
\left[
\dfrac12
\left(
a'(\underline v_k^\varepsilon)
\cdot
\nabla_y\widehat\varphi_k^\varepsilon
\right)
(\mathcal V_k^\varepsilon)^2
\right]
=
0,
\\[2mm]
\mathcal V_k^\varepsilon(0,y,\theta)
=
\mathcal V_{k,0}^\varepsilon(y,\theta).
\end{cases}
\end{equation}

\medskip
\noindent
The presence of the Burgers term in \eqref{eq:local-profile-equation} is
important here. As soon as the periodic profile is nonconstant in the fast
variable and the corresponding quadratic coefficient does not vanish,
shocks are unavoidable in finite time. Indeed, a nonconstant smooth periodic
profile necessarily contains both increasing and decreasing regions. Since
the quadratic Burgers coefficient does not vanish, in one of these regions
neighboring characteristics move towards each other and intersect in finite
time, which produces a shock. We shall check below that, under the
transversality assumption, the quadratic coefficient is uniformly bounded
away from zero.

\medskip
\noindent
Consequently, for fixed $\varepsilon>0$ and $k$, after shock formation
the profile $V_k^\varepsilon$ has to be understood as a $BV$ entropy
solution rather than as a classical one. In particular, the expression
\[
V_k^\varepsilon\left(
    \tau,y,
    \frac{\widehat{\phi}_k^\varepsilon(\tau,y)}{\varepsilon}
\right)
\]
appearing in \eqref{eq:local-WNGO-expansion} is no longer an ordinary
pointwise substitution. Following \cite{Ch92,Ch97,V}, it is understood as the trace
of the $BV$ profile $V_k^\varepsilon$ on the smooth hypersurface
\[
\Sigma_k^\varepsilon
:=
\left\{
(\tau,y,\theta)\,;\,
\theta=
\frac{\widehat{\phi}_k^\varepsilon(\tau,y)}{\varepsilon}
\right\}.
\]
Before shock formation, this trace coincides with the usual pointwise
evaluation. The precise construction of the $BV$ representative and the
justification of this trace after shock formation are given in \cite{Ch92,Ch97}.

\medskip
\noindent
We next check that the Burgers coefficient appearing in the profile equation
is uniformly nonzero. By the transversality assumption and its propagation
given by Lemma~\ref{Lemma propagation}, there exists $c_0>0$ such that
\begin{equation}\label{eq:propagated-transversality}
\left|
a'(\underline u(t,x))\cdot\nabla_x\varphi(t,x)
\right|
\geq c_0
\end{equation}
for every $(t,x)\in\Omega$.

\medskip
\noindent
In the localized variables, the Burgers coefficient is
$b_k^\varepsilon(\tau,y):=
\frac12
a'(\underline v_k^\varepsilon(\tau,y))
\cdot
\nabla_y\widehat\varphi_k^\varepsilon(\tau,y)$.
Since
$\underline v_k^\varepsilon(\tau,y)
=
\underline u(\varepsilon\tau,\varepsilon k+\varepsilon y)$,
and since the local phase $\widehat\varphi_k^\varepsilon$ is the phase
generated by the eikonal equation in the localized variables, the preceding
lower bound gives, uniformly in $k$ and $\varepsilon$,
\begin{equation}\label{eq:local-burgers-coefficient-nonzero}
|b_k^\varepsilon(\tau,y)|
\geq c
\end{equation}
for some $c>0$, for $0\leq \tau\leq\mathcal I$ and $y\in B(0,1]$,
provided $\varepsilon$ is sufficiently small.

\paragraph{Oleinik estimate.}\label{paraOleinik}

The lower bound \eqref{eq:local-burgers-coefficient-nonzero} ensures that the coefficient of the
Burgers-type nonlinear term remains uniformly away from zero. This allows
us, after suitable changes of variables, to reduce the profile equation to
the standard one-dimensional Burgers equation in the periodic variable
$\theta$. The classical Oleinik estimate \cite{lax2005hyperbolic,Ol} can then be applied, yielding the
decay estimate stated below.

\begin{lem}[Oleinik estimate in the fast time variable]
\label{lem:oleinik-local-profile}
Let $\mathcal V_k^\varepsilon$ be the localized profile given by
\eqref{eq:local-profile-equation}. Then, for every $\tau>0$, there exists
a constant $C>0$, independent of $k$ and $\varepsilon$, such that
\[
\int_{B(0,1]}
TV_\theta\left(\mathcal V_k^\varepsilon(\tau,y,\cdot)\right)\,dy
\leq
\frac{C}{\tau}.
\]
\end{lem}

\begin{Proof}
\noindent
We introduce the characteristic flow associated with
$a(\underline v_k^\varepsilon)$. For $z\in B(0,1]$, let
$Y_k^\varepsilon(\tau;z)$ solve
\begin{equation}\label{eq:local-profile-flow}
\begin{cases}
\dfrac{d}{d\tau}Y_k^\varepsilon(\tau;z)
=
a\left(
\underline v_k^\varepsilon
\left(\tau,Y_k^\varepsilon(\tau;z)\right)
\right),
\\[2mm]
Y_k^\varepsilon(0;z)=z.
\end{cases}
\end{equation}
For every $\tau$, the map $z\mapsto Y_k^\varepsilon(\tau;z)$ is a smooth
diffeomorphism onto its image. Denote its Jacobian by
$J_k^\varepsilon:=\det D_zY_k^\varepsilon$. Then
$\partial_\tau J_k^\varepsilon
=\widetilde\alpha_k^\varepsilon J_k^\varepsilon$, where
\[
\widetilde\alpha_k^\varepsilon(\tau,z)
:=
\alpha_k^\varepsilon
\left(\tau,Y_k^\varepsilon(\tau;z)\right),
\qquad
\alpha_k^\varepsilon
:=
\operatorname{div}_y\left(a(\underline v_k^\varepsilon)\right).
\]
We also set
\[
\widetilde{\mathcal V}_k^\varepsilon(\tau,z,\theta)
:=
\mathcal V_k^\varepsilon
\left(\tau,Y_k^\varepsilon(\tau;z),\theta\right),
\qquad
\widetilde b_k^\varepsilon(\tau,z)
:=
b_k^\varepsilon
\left(\tau,Y_k^\varepsilon(\tau;z)\right).
\]

\medskip
\noindent
We now pass to Lagrangian variables in the weak formulation. Let
$\psi\in C_c^1([0,+\infty)\times B(0,1]\times\mathbb T)$ be periodic in
$\theta$, and define
$\phi(\tau,y,\theta)
:=\psi(\tau,(Y_k^\varepsilon(\tau;\cdot))^{-1}(y),\theta)$.
Along the flow,
\[
\left(
\partial_\tau
+
a(\underline v_k^\varepsilon)\cdot\nabla_y
\right)\phi
=
\partial_\tau\psi,
\qquad
\partial_\theta\phi=\partial_\theta\psi.
\]
Using these identities in the weak formulation of
\eqref{eq:local-profile-equation}, together with
$y=Y_k^\varepsilon(\tau;z)$ and
$dy=J_k^\varepsilon(\tau,z)\,dz$, gives
\begin{equation}\label{eq:local-profile-lagrangian-conservative}
\partial_\tau
\left(
J_k^\varepsilon
\widetilde{\mathcal V}_k^\varepsilon
\right)
+
\partial_\theta
\left(
J_k^\varepsilon
\widetilde b_k^\varepsilon
\left(
\widetilde{\mathcal V}_k^\varepsilon
\right)^2
\right)
=
0
\end{equation}
in the sense of distributions. Since $J_k^\varepsilon$ is positive,
independent of $\theta$, and satisfies
$\partial_\tau J_k^\varepsilon
=\widetilde\alpha_k^\varepsilon J_k^\varepsilon$, this is equivalent to
\begin{equation}\label{eq:local-profile-lagrangian}
\partial_\tau
\widetilde{\mathcal V}_k^\varepsilon
+
\widetilde\alpha_k^\varepsilon
\widetilde{\mathcal V}_k^\varepsilon
+
\partial_\theta
\left(
\widetilde b_k^\varepsilon
\left(
\widetilde{\mathcal V}_k^\varepsilon
\right)^2
\right)
=
0.
\end{equation}

\medskip
\noindent
We next eliminate the linear term. Set
$\Gamma_k^\varepsilon(\tau,z)
:=\int_0^\tau\widetilde\alpha_k^\varepsilon(s,z)\,ds$ and
$W_k^\varepsilon
:=e^{\Gamma_k^\varepsilon}\widetilde{\mathcal V}_k^\varepsilon$.
Since $\Gamma_k^\varepsilon$ is independent of $\theta$,
\eqref{eq:local-profile-lagrangian} yields
\begin{equation}\label{eq:local-profile-weighted}
\partial_\tau W_k^\varepsilon
+
\partial_\theta
\left(
\beta_k^\varepsilon(\tau,z)
\left(
W_k^\varepsilon
\right)^2
\right)
=
0,
\end{equation}
where
$\beta_k^\varepsilon
:=\widetilde b_k^\varepsilon e^{-\Gamma_k^\varepsilon}$.
Since $\widetilde b_k^\varepsilon$ is uniformly bounded away from zero and
$\Gamma_k^\varepsilon$ is uniformly bounded on
$0\leq\tau\leq\mathcal I$, there exists $\delta>0$, independent of $k$
and $\varepsilon$, such that
\begin{equation}\label{eq:beta-lower-bound}
\left|
\beta_k^\varepsilon(\tau,z)
\right|
\geq
\delta,
\qquad
0\leq\tau\leq\mathcal I.
\end{equation}

\medskip
\noindent
We first treat the case
$\beta_k^\varepsilon(\tau,z)\geq\delta>0$; the case
$\beta_k^\varepsilon(\tau,z)\leq-\delta$ follows in the same way after
replacing $\theta$ by $-\theta$. For each fixed $z$, define
$\sigma(\tau,z):=\int_0^\tau\beta_k^\varepsilon(s,z)\,ds$.
Then $\tau\mapsto\sigma(\tau,z)$ is a $C^1$-diffeomorphism onto its image and
\[
\delta\tau\leq\sigma(\tau,z)\leq C\tau
\]
for some $C>0$ independent of $k$ and $\varepsilon$.

\medskip
\noindent
Let
$Z_k^\varepsilon(\sigma,z,\theta)
:=W_k^\varepsilon(\tau(\sigma,z),z,\theta)$.
Since
$\partial_\tau W_k^\varepsilon
=\beta_k^\varepsilon\partial_\sigma Z_k^\varepsilon$,
equation~\eqref{eq:local-profile-weighted} reduces to
\begin{equation}\label{eq:local-burgers-reduced}
\partial_\sigma Z_k^\varepsilon
+
\partial_\theta
\left(
\left(Z_k^\varepsilon\right)^2
\right)
=
0.
\end{equation}
Thus, for every fixed $z$, $Z_k^\varepsilon$ solves the one-dimensional
Burgers equation in $(\sigma,\theta)$.

\medskip
\noindent
By the Oleinik estimate for one-dimensional scalar conservation laws with
uniformly convex flux,
\[
TV_\theta
\left(
Z_k^\varepsilon(\sigma,z,\cdot)
\right)
\leq
\frac{C}{\sigma}.
\]
Since $\sigma(\tau,z)\geq\delta\tau$, we obtain
\[
TV_\theta
\left(
W_k^\varepsilon(\tau,z,\cdot)
\right)
\leq
\frac{C}{\tau}.
\]
Finally,
$\widetilde{\mathcal V}_k^\varepsilon
=e^{-\Gamma_k^\varepsilon}W_k^\varepsilon$, and
$\Gamma_k^\varepsilon$ is independent of $\theta$ and uniformly bounded.
Hence
\[
TV_\theta
\left(
\widetilde{\mathcal V}_k^\varepsilon(\tau,z,\cdot)
\right)
\leq
\frac{C}{\tau}.
\]
Returning to $y=Y_k^\varepsilon(\tau;z)$ and using the uniform upper and
lower bounds on $J_k^\varepsilon$, we finally obtain
\begin{equation}\label{eq:local-oleinik-final}
\int_{B(0,1]}
TV_\theta
\left(
\mathcal V_k^\varepsilon(\tau,y,\cdot)
\right)
\,dy
\leq
\frac{C}{\tau},
\qquad
0<\tau\leq\mathcal I.
\end{equation}
\end{Proof}

\subsubsection{Comparison with the averaged profile}
\label{compaave}

We now use the Oleinik estimate to show that the oscillating profile is close
to its average with respect to the fast variable. We define
\[
\overline{\mathcal V}_k^\varepsilon(\tau,y)
:=
\int_{\mathbb T}
\mathcal V_k^\varepsilon(\tau,y,\theta)\,d\theta .
\]
Since $\mathcal V_k^\varepsilon$ is periodic with respect to $\theta$, for
every $\theta\in\mathbb T$,
\[
\left|
\mathcal V_k^\varepsilon(\tau,y,\theta)
-
\overline{\mathcal V}_k^\varepsilon(\tau,y)
\right|
\leq
TV_\theta
\left(
\mathcal V_k^\varepsilon(\tau,y,\cdot)
\right).
\]
Therefore, by Lemma~\ref{lem:oleinik-local-profile}, taking
$\tau=\mathcal I$ gives
\begin{equation}\label{eq:profile-close-average}
\int_{B(0,1]}
\left|
\mathcal V_k^\varepsilon
\left(
\mathcal I,y,
\frac{\widehat\varphi_k^\varepsilon(\mathcal I,y)}{\varepsilon}
\right)
-
\overline{\mathcal V}_k^\varepsilon(\mathcal I,y)
\right|
\,dy
\leq
\frac{C}{\mathcal I}.
\end{equation}

\medskip
\noindent
We now identify this averaged profile with the rescaled averaged correction.
Integrating \eqref{eq:local-profile-equation} with respect to $\theta$, the
$\partial_\theta$-term vanishes by periodicity, so
\[
\partial_\tau \overline{\mathcal V}_k^\varepsilon
+
\operatorname{div}_y
\left(
a(\underline v_k^\varepsilon)
\overline{\mathcal V}_k^\varepsilon
\right)
=
0.
\]
Moreover, from the definition of $\mathcal V_{k,0}^\varepsilon$, the phase
shift is only a translation in the periodic variable, and therefore
\[
\overline{\mathcal V}_k^\varepsilon(0,y)
=
\int_{\mathbb T}
U_0(\varepsilon k+\varepsilon y,\theta)\,d\theta
=
\overline U_0(\varepsilon k+\varepsilon y).
\]
On the other hand,
$\overline U_k^\varepsilon(\tau,y)
:=\overline U(\varepsilon\tau,\varepsilon k+\varepsilon y)$
solves the same linear equation with the same initial datum. Hence, by
uniqueness,
\begin{equation}\label{eq:average-profile-identification}
\overline{\mathcal V}_k^\varepsilon
=
\overline U_k^\varepsilon.
\end{equation}

\medskip
\noindent
We now combine this estimate with the WNGO approximation. By
\cref{eq:local-WNGO-expansion} at $\tau=\mathcal I$, together with
\eqref{eq:profile-close-average} and
\eqref{eq:average-profile-identification}, the triangle inequality gives
\[
\begin{aligned}
&
\left\|
v_k^\varepsilon(\mathcal I,\cdot)
-
\underline v_k^\varepsilon(\mathcal I,\cdot)
-
\varepsilon
\overline U_k^\varepsilon(\mathcal I,\cdot)
\right\|_{L^1(B(0,1])}
\\
&\qquad
\leq
C_{\mathcal I}\varepsilon^2
+
\varepsilon
\int_{B(0,1]}
\left|
\mathcal V_k^\varepsilon
\left(
\mathcal I,y,
\frac{\widehat\varphi_k^\varepsilon(\mathcal I,y)}{\varepsilon}
\right)
-
\overline{\mathcal V}_k^\varepsilon(\mathcal I,y)
\right|
\,dy
\\
&\qquad
\leq
C_{\mathcal I}\varepsilon^2
+
\frac{C\varepsilon}{\mathcal I}.
\end{aligned}
\]
Thus
\begin{equation}\label{eq:local-cancellation-estimate-before-choice}
\left\|
v_k^\varepsilon(\mathcal I,\cdot)
-
\underline v_k^\varepsilon(\mathcal I,\cdot)
-
\varepsilon
\overline U_k^\varepsilon(\mathcal I,\cdot)
\right\|_{L^1(B(0,1])}
\leq
\frac{C\varepsilon}{\mathcal I}
+
C_{\mathcal I}\varepsilon^2 .
\end{equation}

\medskip
\noindent
Since $\mathcal I$ is fixed, we choose $\varepsilon_{\mathcal I}>0$ small
enough so that
$C_{\mathcal I}\varepsilon^2\leq C\varepsilon/\mathcal I$
for every $0<\varepsilon\leq\varepsilon_{\mathcal I}$.
Then, enlarging $C$ if necessary,
\[
\left\|
v_k^\varepsilon(\mathcal I,\cdot)
-
\underline v_k^\varepsilon(\mathcal I,\cdot)
-
\varepsilon
\overline U_k^\varepsilon(\mathcal I,\cdot)
\right\|_{L^1(B(0,1])}
\leq
C\frac{\varepsilon}{\mathcal I}.
\]
This proves Proposition~\ref{prop:local-cancellation}.

\section{Amplitude $\varepsilon$ with frequency $\varepsilon^{-p}$ when $p>2$}
\label{sec:3}

\medskip
\noindent
This section is organized into five subsections.
In Subsection~\ref{subsec:compatibility-algebra-rank}, we study more closely the
Lie algebra $ \mathcal A_p $; we explain the propagation of the
constant-rank condition; and we discuss a few situations illustrating its
geometrical meaning. In Subsection~\ref{subsec:sufficient-condition}, we
prove the sufficient part of 
Theorem~\ref{thm:complete-WKB-criterion} by constructing a compatible
non-trivial phase and the complete WKB hierarchy. In
Subsection~\ref{subsec:necessary-condition}, we prove the converse
implication by recovering the compatibility conditions and by excluding
the full-rank case. Subsection~\ref{subsec:explicit-geometric-criteria}
is devoted to particular situations in which the rank 
condition can be described more explicitly. Finally, in
Subsection~\ref{Cancellation-p>2}, we work 
under Assumption~\ref{ass-fullyn}; we study the cancellation
regime; and we justify 
Theorem~\ref{thm-cancel3}.


\subsection{Lie algebra and rank properties}
\label{subsec:compatibility-algebra-rank}

\noindent
We examine more closely the Lie algebra generated by the $ L_k $  
and the behavior of its rank. We first show that the constant-rank
condition imposed at $ t = 0 $ is preserved by the evolution of
the smooth ground state. We then present two  examples illustrating
the role of this assumption. 


\subsubsection{Propagation of the rank condition}
\label{subsubsec:propagation-rank-condition}

\noindent
Assumption~\ref{ass:constant-rank-Np} is  propagated by the characteristic
flow of the smooth ground state. Indeed:

\begin{lem}[Propagation of the rank condition]
\label{lem:propagation-rank-Ap} Under Assumption \ref{ass:constant-rank-Np}, we have
\begin{equation}\label{consttttN}
N_p(t,x)=N_p,
\qquad \forall \, 
(t,x)\in\Omega \, .
\end{equation}
\end{lem}

\begin{Proof}
Let $ T :=\partial_t+a(\underline u)\cdot\nabla_x $
be the transport operator associated with the smooth ground state, and
let $\Phi_t$ be its characteristic flow. For a spatial vector field
$X(t,\cdot)$, we denote by $\Phi_t^*X(t)$ its pullback by the flow,
namely
\[
\bigl(\Phi_t^*X(t)\bigr)(y)
:=
D\Phi_t(y)^{-1}
X\bigl(t,\Phi_t(y)\bigr).
\]
Since $ T \underline u=0$, a direct computation gives
\begin{equation}\label{eq:D-Lk-commutator}
[T,L_k] = - ( L_k \underline u ) \, L_1 ,
\qquad
k=1,\ldots,p-1.
\end{equation}
We now pull the vector fields back by the flow and set $ Y_k(t,y):= (\Phi_t^*L_k)(t,y) $.
Using the standard differentiation formula for time-dependent vector
fields along a flow, \eqref{eq:D-Lk-commutator} yields
\begin{equation}\label{eq:Yk-evolution}
\frac{d}{dt}Y_k(t,y)
=\frac{d}{dt} (\Phi_t^*L_k)(t,y) = \Phi_t^* \, [ T,L_k] (t,y)
 =-\gamma_k(t,y) \, Y_1(t,y),
\end{equation}
where 
\[
\gamma_k(t,y)
= \Phi_t^* \, (L_k\underline u) (t,y) = 
\bigl(L_k\underline u\bigr)\bigl(t,\Phi_t(y)\bigr).
\]
For $k=1$, this is a scalar linear equation starting from $ Y_1(0,y) = L_1(0,y) $, hence
\[
Y_1(t,y)=\alpha(t,y)L_1(0,y),
\qquad
\alpha(t,y)\neq0.
\]
For $k\geq2$, this relation together with the integration of \eqref{eq:Yk-evolution} furnishes
\[
Y_k(t,y)
=
L_k(0,y)+\beta_k(t,y)L_1(0,y) \, , \qquad \beta_k(t,y) := - \int_0^t (\gamma_k \, \alpha)(s,y) \ ds \, .
\]
It follows that the families
$Y_1(t,y),\ldots,Y_{p-1}(t,y)$ and
$L_1(0,y),\ldots,L_{p-1}(0,y)$ generate the same Lie algebra.
Since pullback by a diffeomorphism commutes with Lie brackets, their
Lie closures coincide, so that $ \Phi_t^*\mathcal A_p(t)=\mathcal A_p(0) $, and therefore $ N_p\bigl(t,\Phi_t(y)\bigr)=N_p(0,y) $. Since by
assumption, the initial rank is constant and equal to $N_p$ on $\Omega_0$, the same applies at time $ t $.
\end{Proof}

\subsubsection{Two representative scenarios}
\label{subsubsec:comments-rank-assumption}

\noindent This paragraph provides an opportunity to revisit two points mentioned earlier. The former addresses the inclusion $ \mathcal D_p \subset \mathcal A_p $, whereas the latter pertains to Assumption \ref{ass:constant-rank-Np}.

\smallskip

\noindent $ \bullet $ {\it Distinction between $ \mathcal D_p $ and $ \mathcal A_p $.} As said before, the distributions $ \mathcal D_p $ and $ \mathcal A_p $ coincide when the ground state $ \underline u $ is constant. However, for a nonconstant ground state $ \underline u $, the situation may be quite different.
As a matter of fact, the Lie brackets usually create directions
which are not contained inside $\mathcal D_p$. The Lie algebra $ \mathcal A_p $ may have full rank $ d $ even
when $p \leq d$. To understand why, take for instance $p=d=3$, and consider a flux function $ F(u) $ such that
\[
F'(u) = a(u)=
\left(
u,\frac{u^2}{2},\frac{u^3}{6}
\right),
\qquad
\underline u_0(x)=x_1.
\]
At the origin,
\[
L_1(0)=\partial_{x_1},
\qquad
L_2(0)=\partial_{x_2},
\qquad
[L_1,L_2](0)=\partial_{x_3},
\]
and therefore $ \text{dim} \, \mathcal D_3 (0,0) = 2 < 3 $ while $
N_3(0,0)=3=d $. This example confirms that the restriction $p\leq d$, which is sufficient
to leave transverse directions in the constant-background setting, is
no longer sufficient for a nonconstant ground state.

\smallskip

\noindent $ \bullet $ {\it The role of the rank condition.} Assumption \ref{ass:constant-rank-Np}, which is necessary in order to apply Frobenius Theorem, is not always satisfied. Again with $p=d=3$, consider a flux function $ F(u) $ such that
\[ F'(u) = 
a(u)=
\left(
u,\frac{u^3}{3},0
\right),
\qquad
\underline u_0(x)=x_2.
\]
Then
\[
L_1
=
\partial_{x_1} + x_2^2 \, \partial_{x_2},
\qquad
L_2
=
2 \, x_2 \, \partial_{x_2} . 
\]
Observe that 
\[ [L_1,L_2] = - 2 \, x_2^2 \, \partial_{x_2} \, , \qquad \bigl \lbrack L_1, [L_1,L_2] \bigr \rbrack = 0 \, , \qquad \bigl \lbrack L_2, [L_1,L_2] \bigr \rbrack = 2 \, [L_1,L_2] \, . \] 
It follows that $ \text{Lie} (L_1,L_2) $ is generated by $ \partial_{x_1} $ when $ x_2 = 0 $, so  that $ N_3(x_1,0)=1 $; and it is spanned by $ \partial_{x_1} $ and $ \partial_{x_2} $ when $ x_2 \not = 0 $, in which case $ N_3(x_1,x_2)=2 $.
Observe that, in any open neighbourhood of the origin, the rank cannot be  constant. Still, there is a compatible phase, namely $ \phi(x) = x_3 $. Yet, there are situations for which the rank is not locally constant (and therefore everywhere $ <d $), while there exists no non-trivial phase. 

\medskip
\noindent
Such variable-rank situations are covered by the Stefan--Sussmann theory of
singular distributions and orbits of families of smooth vector fields; see,
for instance, \cite{Sussmann73}. We do not pursue this more general setting
here, since our analysis relies on the constant-rank
Assumption~\ref{ass:constant-rank-Np}.


\subsection{The sufficient condition ($N_p<d$)}
\label{subsec:sufficient-condition}

\noindent
\noindent
We prove first the sufficient part in the implication of
Theorem~\ref{thm:complete-WKB-criterion}. Throughout this subsection,
we assume that $ N_p<d $.
The constant-rank condition and the Frobenius theorem allow us to
construct a non-trivial initial phase compatible with the algebra
$\mathcal A_p(0)$. We then propagate the corresponding compatibility
conditions by the eikonal equation, construct successively the WKB
profiles, and finally establish the uniform estimates on the
remainder required in \eqref{eq:complete-WKB-remainder}.

\subsubsection{Construction and propagation of compatible phases}
\label{subsubsec:construction-compatible-phases}

\noindent
The goal here is to construct a non-trivial phase compatible with
$\mathcal A_p$. The initial phase is obtained from the Frobenius
theorem (Paragraph \ref{par:construction-initial-compatible-phase}), and its compatibility conditions are then propagated by the
eikonal evolution (Paragraph \ref{par:propagation-compatibility-conditions}).


\paragraph{Construction of the initial compatible phase.}
\label{par:construction-initial-compatible-phase}

\noindent
Since $\mathcal A_p(0,\cdot)$ has constant rank $N_p<d$ and is
involutive, Frobenius theorem gives, after possibly reducing the
neighborhood of the origin, local coordinates $ (y,z)\in\mathbb R^{N_p}\times\mathbb R^{d-N_p} $ such that
\[
\mathcal A_p(0)
=
\operatorname{span}
\left\{
\frac{\partial}{\partial y_1},
\ldots,
\frac{\partial}{\partial y_{N_p}}
\right\}.
\]
We can choose
\begin{equation}\label{eq:phi0-adapted-coordinates}
\phi_0(y,z)=\Phi(z),
\end{equation}
with $\nabla_z\Phi\neq0$. Hence
\begin{equation}\label{eq:Ap0-annihilates-phi0}
T\phi_0=0,
\qquad
\forall \, T\in\mathcal A_p(0),
\end{equation}
while $\nabla_x\phi_0\neq0$. In particular,
\begin{equation}\label{eq:initial-compatibility-p}
L_k(0)\phi_0=0,
\qquad
k=1,\ldots,p-1.
\end{equation}


\paragraph{Propagation of the compatibility conditions.}
\label{par:propagation-compatibility-conditions}

\noindent We claim that a phase $\phi$ 
which is compatible at time $ t = 0 $ remains compatible for $ 0 < t $.

\begin{lem}
\label{Lemma propagation Lie algebra}
Let $\phi$ be the solution on $ \Omega $ of the eikonal equation
\eqref{eq:eikonal} with initial condition
$\phi|_{t=0}=\phi_0$ satisfying \eqref{eq:initial-compatibility-p}. Then, 
for every $t\in[0,T]$, we have
\begin{equation}\label{eq:propagated-generators}
L_k(t)\phi(t,\cdot)=0,
\qquad
k=1,\ldots,p-1.
\end{equation}
Moreover,
\begin{equation}\label{TLielie}
T\phi(t,\cdot)=0,
\qquad
\forall \, T \in
\operatorname{Lie}\bigl(L_1(t),\ldots,L_{p-1}(t)\bigr).
\end{equation}
\end{lem}

\begin{Proof} 
By construction of the initial phase, \eqref{eq:initial-compatibility-p} holds.
For $k=1$, part~{\rm (i)} of Lemma~\ref{Lemma propagation}, together with
\eqref{eq:initial-compatibility-p}, gives rise to
\[
a'(\underline u(t,x))\cdot\nabla_x\phi(t,x)=0,
\qquad
(t,x)\in\Omega.
\]
Hence $ L_1(t)\phi(t,\cdot)=0 $. For $k=2,\ldots,p-1$, part~{\rm (ii)} of
Lemma~\ref{Lemma propagation} with $ G=a^{(k)} $. Since
\[
a^{(k)}(\underline u_0)\cdot\nabla_x\phi_0
=
L_k(0)\phi_0
=
0
\]
by \eqref{eq:initial-compatibility-p}, we obtain
\[
a^{(k)}(\underline u(t,x))\cdot\nabla_x\phi(t,x)=0,
\qquad
(t,x)\in\Omega.
\]
Therefore, \eqref{eq:propagated-generators} holds.

\medskip
\noindent
It remains to extend this property from the generators to the whole
compatibility algebra. Fix $t\in[0,T]$. If two smooth vector fields
$X$ and $Y$ satisfy
\[
X\phi(t,\cdot)=0,
\qquad
Y\phi(t,\cdot)=0,
\]
then every smooth linear combination of $X$ and $Y$ also annihilates
$\phi(t,\cdot)$. Moreover,
\[
[X,Y]\phi(t,\cdot)
=
X\bigl(Y\phi(t,\cdot)\bigr)
-
Y\bigl(X\phi(t,\cdot)\bigr)
=
0.
\]
Thus the property of annihilating $\phi(t,\cdot)$ is stable under smooth
linear combinations and Lie brackets. Starting from
\eqref{eq:propagated-generators}, we recover \eqref{TLielie}.
\end{Proof}


\subsubsection{Construction of the complete WKB hierarchy}
\label{subsubsec:complete-WKB-hierarchy}

\noindent
With the compatible phase constructed above, we now derive the complete
WKB hierarchy. The compatibility conditions remove the singular terms
generated by the fast variable and allow the profiles to be constructed
successively. The principal profile satisfies a nonlinear equation,
while the higher-order profiles satisfy linear equations.


\paragraph{The complete WKB ansatz.}
\label{par:complete-WKB-ansatz}

\noindent We start with a formal WKB calculus \cite{Rau12}, assuming that all profiles are smooth. Fix $K\in\mathbb N^*$. We seek an approximate solution in the form
\begin{equation}\label{eq:truncated-WKB}
v_K^\varepsilon(t,x)
=
\underline u(t,x)
+
\sum_{j=1}^{K}
\varepsilon^j \, 
U^j\!\left(
t,x,\frac{\phi(t,x)}{\varepsilon^p}
\right).
\end{equation}
We set $ \theta=\phi(t,x)/\varepsilon^p $,
and we introduce the truncated profile
\[
W_K^\varepsilon(t,x,\theta)
:=
\sum_{j=1}^{K}\varepsilon^j \, U^j(t,x,\theta).
\]
Thus,
\[
v_K^\varepsilon(t,x)
=
\underline u(t,x)
+
W_K^\varepsilon
\left(
t,x,\frac{\phi(t,x)}{\varepsilon^p}
\right).
\]
Since the profiles are smooth and bounded on the compact
domain $ \Omega $, we have
$ W_K^\varepsilon=O(\varepsilon) $ 
uniformly with respect to $(t,x,\theta)$.
We next expand the flux around $ \underline u $. Differentiation through
the variable $\theta=\phi/\varepsilon^p$ produces a factor
$\varepsilon^{-p}$. Hence the Taylor expansion has to be carried
sufficiently far so that its remainder is still of order at least
$\varepsilon^{K+1}$ after differentiation. We therefore expand $F$ up
to degree $K+p$. Taylor's formula gives
\begin{equation}\label{eq:flux-Taylor-WKB}
F(\underline u+W_K^\varepsilon)
=
F(\underline u)
+
\sum_{m=1}^{K+p}
\frac{1}{m!} \, 
a^{(m-1)}(\underline u) \, 
\bigl(W_K^\varepsilon\bigr)^m
+
R_{K+p+1}^\varepsilon,
\end{equation}
where
\begin{equation}\label{eq:flux-Taylor-remainder-WKB}
R_{K+p+1}^\varepsilon
=
O\!\left(
|W_K^\varepsilon|^{K+p+1}
\right)
=
O(\varepsilon^{K+p+1}) \, .
\end{equation}
Since the profiles and their derivatives are uniformly bounded, the Taylor remainder is such that
\begin{equation}\label{taylorremainder}
\partial_{x_j} R_{K+p+1}^\varepsilon
=
O(\varepsilon^{K+p+1})
+
\varepsilon^{-p} \, O(\varepsilon^{K+p+1})
=
O(\varepsilon^{K+1}) \, , \qquad \forall \, j \in \{1,\cdots , d \} \, .
\end{equation}
We now substitute \eqref{eq:truncated-WKB} into the conservation law.
For a profile $V=V(t,x,\theta)$ evaluated at
$\theta=\phi/\varepsilon^p$, the chain rule gives
\[
\partial_t
V\!\left(t,x,\frac{\phi}{\varepsilon^p}\right)
=
\left(
\partial_tV
+
\varepsilon^{-p}
\partial_t\phi\,\partial_\theta V
\right)_{\theta=\phi/\varepsilon^p},
\]
and
\[
\nabla_x
V\!\left(t,x,\frac{\phi}{\varepsilon^p}\right)
=
\left(
\nabla_xV
+
\varepsilon^{-p} \, 
\nabla_x\phi \,\partial_\theta V
\right)_{\theta=\phi/\varepsilon^p}.
\]
Applying these identities term by term to
\eqref{eq:flux-Taylor-WKB}, we obtain
\begin{align}
\partial_t v_K^\varepsilon {}&
+
\operatorname{div}_xF(v_K^\varepsilon)
=
\partial_t\underline u
+
\operatorname{div}_xF(\underline u)
\nonumber\\
&+
\partial_tW_K^\varepsilon
+
\operatorname{div}_x
\left(
a(\underline u)W_K^\varepsilon
\right)
\nonumber\\
&+
\sum_{m=2}^{K+p}
\frac{1}{m!}
\operatorname{div}_x
\left(
a^{(m-1)}(\underline u)
(W_K^\varepsilon)^m
\right)
\nonumber\\
&+
\varepsilon^{-p}
\left(
\partial_t\phi
+
a(\underline u)\cdot\nabla_x\phi
\right)
\partial_\theta W_K^\varepsilon
\nonumber\\
&+
\varepsilon^{-p}
\sum_{m=2}^{K+p}
\frac{1}{m!}
\left(
a^{(m-1)}(\underline u)\cdot\nabla_x\phi
\right)
\partial_\theta
\left(
(W_K^\varepsilon)^m
\right) +
O(\varepsilon^{K+1}),
\label{eq:full-WKB-expansion}
\end{align}
where all profile terms on the right-hand side are evaluated at
$\theta=\phi/\varepsilon^p$.
The first line of \eqref{eq:full-WKB-expansion} disappears because
$\underline u$ solves \cref{eq:claw}, and the fourth line vanishes
because $\phi$ solves the eikonal equation \eqref{eq:eikonal}.
Moreover, by the propagated compatibility conditions
\eqref{eq:propagated-generators}, the terms in the last sum
corresponding to $m=2,\ldots,p$ vanish. Consequently,
\eqref{eq:full-WKB-expansion} reduces to
\begin{align}
\partial_t v_K^\varepsilon {}&
+
\operatorname{div}_xF(v_K^\varepsilon)
=
\partial_tW_K^\varepsilon
+
\operatorname{div}_x
\left(
a(\underline u)W_K^\varepsilon
\right)
\nonumber\\
&+
\sum_{m=2}^{K+p}
\frac{1}{m!}
\operatorname{div}_x
\left(
a^{(m-1)}(\underline u)
(W_K^\varepsilon)^m
\right)
\nonumber\\
&+
\varepsilon^{-p}
\sum_{m=p+1}^{K+p}
\frac{1}{m!}
\left(
a^{(m-1)}(\underline u)\cdot\nabla_x\phi
\right)
\partial_\theta
\left(
(W_K^\varepsilon)^m
\right)+
O(\varepsilon^{K+1}).
\label{eq:reduced-WKB-expansion}
\end{align}

\noindent
The first term in the last sum corresponds to $m=p+1$ and is
\begin{equation}\label{eq:first-surviving-fast-term}
\frac{\varepsilon^{-p}}{(p+1)!}
\bigl(
a^{(p)}(\underline u)\cdot\nabla_x\phi
\bigr)
\partial_\theta
\left(
(W_K^\varepsilon)^{p+1}
\right).
\end{equation}
Moreover,
\begin{equation}\label{eq:WK-leading-order}
W_K^\varepsilon
=
\varepsilon \, U^1+O(\varepsilon^2).
\end{equation}
Hence
\[
(W_K^\varepsilon)^{p+1}
=
\varepsilon^{p+1} \, (U^1)^{p+1}
+
O(\varepsilon^{p+2}),
\]
and therefore
\[
\varepsilon^{-p} \, 
(W_K^\varepsilon)^{p+1}
=
\varepsilon \, (U^1)^{p+1}
+
O(\varepsilon^2).
\]
Thus all the terms remaining in \eqref{eq:reduced-WKB-expansion} have
non-negative powers of $\varepsilon$, and the first non-zero order is
$\varepsilon$. We now collect the coefficient of $\varepsilon$ in
\eqref{eq:reduced-WKB-expansion}. By
\eqref{eq:WK-leading-order}, the first two terms contribute according to
\[
\partial_tU^1
+
\operatorname{div}_x
\bigl(
a(\underline u) \, U^1
\bigr).
\]
The terms
\[
\operatorname{div}_x
\left(
a^{(m-1)}(\underline u) \, 
(W_K^\varepsilon)^m
\right),
\qquad m\geq2,
\]
start at order $\varepsilon^2$, since
$ (W_K^\varepsilon)^m=O(\varepsilon^m) $.
In the last sum, only the term $m=p+1$ contributes at order
$\varepsilon$, and its contribution is
\[
\frac{1}{(p+1)!} \, 
\bigl(
a^{(p)}(\underline u)\cdot\nabla_x\phi
\bigr) \, 
\partial_\theta \bigl( (U^1)^{p+1} \bigr) .
\]
To suppress the contribution with $ \varepsilon $ in factor, we have therefore to impose
\begin{equation}\label{eq:U1-equation}
\partial_tU^1
+
\operatorname{div}_x\bigl(a(\underline u) \, U^1\bigr)
+
\frac{1}{(p+1)!} \, 
\bigl(
a^{(p)}(\underline u)\cdot\nabla_x\phi
\bigr) \, 
\partial_\theta \bigl( (U^1)^{p+1} \bigr)
=
0.
\end{equation}
The initial condition is chosen as
\begin{equation}\label{eq:U1-initial}
U^1(0,x,\theta)=U^1_0(x,\theta).
\end{equation}
We next determine successively the higher profiles. Assume that
$U^1,\ldots,U^{n-1}$ have already been constructed, with
$2\leq n\leq K$. We collect the coefficient of $\varepsilon^n$ in
\eqref{eq:reduced-WKB-expansion}.
To organize the discussion, we introduce the notation
\begin{equation}\label{eq:Pmq-def}
P_{m,q}^{K}
:=
\sum_{\substack{
j_1+\cdots+j_m=q\\
1\leq j_\ell\leq K}}
U^{j_1}\cdots U^{j_m},
\end{equation}
so that
\begin{equation}\label{eq:WK-power-expansion}
(W_K^\varepsilon)^m
=
\sum_{q=m}^{mK}
\varepsilon^q \, P_{m,q}^{K}.
\end{equation}
Indeed,
\[
(W_K^\varepsilon)^m
=
\left(
\sum_{j=1}^{K}\varepsilon^j \, U^j
\right)^m
=
\sum_{j_1,\ldots,j_m=1}^{K}
\varepsilon^{j_1+\cdots+j_m} \, 
U^{j_1}\cdots U^{j_m},
\]
and $P_{m,q}^K$ collects precisely the terms for which
$j_1+\cdots+j_m=q$.
The first two terms in \eqref{eq:reduced-WKB-expansion} contribute
\[
\partial_tU^n
+
\operatorname{div}_x
\bigl(
a(\underline u) \, U^n
\bigr)
\]
at order $\varepsilon^n$.
The terms
\[
\operatorname{div}_x
\left(
a^{(m-1)}(\underline u) \, 
(W_K^\varepsilon)^m
\right),
\qquad m\geq2,
\]
contribute at order $\varepsilon^n$ through the coefficients
$P_{m,n}^K$. However, $P_{m,n}^K$ cannot contain $U^n$. Indeed, if one
of the indices in
\[
j_1+\cdots+j_m=n
\]
were equal to $n$, then the remaining $m-1$ indices would have to sum
to zero, which is impossible since each index is at least one.
Therefore $P_{m,n}^K$ depends only on
$U^1,\ldots,U^{n-1}$.

\smallskip
\noindent
Among the terms carrying the factor $\varepsilon^{-p}$, a contribution
of order $\varepsilon^n$ must arise from a term of order
$\varepsilon^{n+p}$ in $(W_K^\varepsilon)^m$. By
\eqref{eq:WK-power-expansion}, the relevant coefficient is therefore
$P_{m,n+p}^K$. Since the smallest power occurring in
$(W_K^\varepsilon)^m$ is $\varepsilon^m$, only the indices $ m=p+1,\ldots,n+p $ 
can contribute at order $\varepsilon^n$.
We first determine which of these contributions may contain the unknown
profile $U^n$. Suppose that one of the indices entering
$P_{m,n+p}^K$ is equal to $n$. The remaining $m-1$ indices must then
satisfy $ j_1+\cdots+j_{m-1}=p $.
Since every index is at least one, this requires $ m-1\leq p $.
On the other hand, $m\geq p+1$. Consequently, the only possible value
is $ m=p+1 $. Thus the only contribution at order $\varepsilon^n$ which can contain
$U^n$ comes from $P_{p+1,n+p}^K$.

\smallskip
\noindent
We now isolate the terms of $P_{p+1,n+p}^K$ containing $U^n$. If one of
the $p+1$ indices is equal to $n$, the remaining $p$ indices must have
sum $p$. Since each of them is at least one, they are necessarily all
equal to one. There are $p+1$ possible positions for the factor $U^n$.
Hence
\begin{equation}\label{eq:P-linear-Un}
P^{K}_{p+1,n+p}
=
(p+1) \, (U^1)^p \, U^n
+
\widetilde P^{K}_{p+1,n+p},
\end{equation}
where $\widetilde P^{K}_{p+1,n+p}$ depends only on
$U^1,\ldots,U^{n-1}$.

\smallskip
\noindent
For $m=p+2,\ldots,n+p$, the coefficients $P_{m,n+p}^K$ may also
contribute to the order-$\varepsilon^n$ equation, but none of them can
contain $U^n$. Indeed, if one index were equal to $n$, the remaining
$m-1\geq p+1$ positive indices would have to sum to $p$, which is
impossible. Therefore all these contributions depend only on the
previously constructed profiles $U^1,\ldots,U^{n-1}$. Consequently, the contribution of order $\varepsilon^n$ which involves $U^n$
reduces to 
\[
\partial_tU^n
+
\operatorname{div}_x
\bigl(
a(\underline u) \, U^n
\bigr)
+
\frac{1}{p!} \, 
\bigl(
a^{(p)}(\underline u)\cdot\nabla_x\phi
\bigr) \, 
\partial_\theta \bigl( (U^1)^p \, U^n \bigr) .
\]
All the remaining terms of order $\varepsilon^n$ depend only on the
already constructed profiles $U^1,\ldots,U^{n-1}$, and we collect them
in a source term $S^n$. Thus
\begin{equation}\label{eq:Un-linear-equation}
\mathcal L_{U^1}U^n
=
S^n(U^1,\ldots,U^{n-1}),
\end{equation}
where
\begin{equation}\label{eq:L-U1}
\mathcal L_{U^1}V
=
\partial_t V
+
\operatorname{div}_x\bigl(a(\underline u) \, V\bigr)
+
\frac{1}{p!} \, \bigl(
a^{(p)}(\underline u)\cdot\nabla_x\phi
\bigr) \, 
\partial_\theta \bigl( (U^1)^p \, V \bigr) .
\end{equation}
Once $U^1$ is known and smooth, the operator $\mathcal L_{U^1}$ is linear with
respect to $V$, and therefore \eqref{eq:Un-linear-equation} determines
$U^n$ from the initial data $ U^n (0,\cdot) = U^n_0 $ successively from the previously constructed profiles. At the end, there remains
\begin{equation}\label{eq:physical-WKB-residual}
\partial_t v_K^\varepsilon
+
\operatorname{div}_xF(v_K^\varepsilon)
=
\varepsilon^{K+1} \, R_K^\varepsilon,
\end{equation}
where $R_K^\varepsilon$ is uniformly bounded and subject to \eqref{taylorremainder}.


\paragraph{Exact $\varepsilon$-dependent profile.}
\label{par:exact-epsilon-profile}

\noindent
We now turn to the comparison between the exact oscillating solution
and the truncated WKB approximation. To this end, we introduce an exact
$\varepsilon$-dependent profile $U^\varepsilon(t,x,\theta)$ through the
reconstruction
\begin{equation}\label{eq:exact-profile-reconstruction}
u^\varepsilon(t,x)
=
\underline u(t,x)
+
\varepsilon
U^\varepsilon\!\left(
t,x,\frac{\phi(t,x)}{\varepsilon^p}
\right),
\end{equation}
which is required to satisfy exactly the scalar conservation law
\begin{equation}\label{eq:exact-SCL}
\partial_tu^\varepsilon
+
\operatorname{div}_xF(u^\varepsilon)
=
0.
\end{equation}
At $t=0$, the reconstruction must agree with the prescribed
oscillatory initial data, so that
\begin{equation}\label{eq:exact-profile-initial}
U^\varepsilon(0,x,\theta)
= \sum_{k=1}^K \varepsilon^{k-1} \, U_0^k \!\left(x,\frac{\phi_0(x)}{\varepsilon^p}\right) .
\end{equation}
On the other hand, the truncated WKB approximation can be written as
\[
v_K^\varepsilon(t,x)
=
\underline u(t,x)
+
\varepsilon
U_{[K]}^\varepsilon\!\left(
t,x,\frac{\phi(t,x)}{\varepsilon^p}
\right),
\]
where we set
\begin{equation}\label{eq:truncated-amplitude-profile}
U_{[K]}^\varepsilon(t,x,\theta)
:=
U^1(t,x,\theta)
+
\varepsilon \, U^2(t,x,\theta)
+\cdots+
\varepsilon^{K-1} \, U^K(t,x,\theta).
\end{equation}
Thus the comparison between the exact solution $u^\varepsilon$ and the
approximate solution $v_K^\varepsilon$ is reduced to a comparison
between the exact profile $U^\varepsilon$ and the truncated profile
$U_{[K]}^\varepsilon$. Indeed,
\begin{equation}\label{eq:physical-profile-difference}
u^\varepsilon(t,x)-v_K^\varepsilon(t,x)
=
\varepsilon
\left[
U^\varepsilon
-
U_{[K]}^\varepsilon
\right]
\left(
t,x,\frac{\phi(t,x)}{\varepsilon^p}
\right).
\end{equation}
The next step is therefore to derive the exact equation satisfied by
$U^\varepsilon$ and to compare it with the equation satisfied, up to a
small residual, by $U_{[K]}^\varepsilon$.


\paragraph{Profile equations.}
\label{par:profile-equations}

\noindent
By the same chain-rule computation as in the WKB expansion, the exact
profile $U^\varepsilon$ introduced in
\eqref{eq:exact-profile-reconstruction} satisfies
\begin{equation}\label{eq:exact-profile-equation}
\partial_tU^\varepsilon
+
\operatorname{div}_y
G^\varepsilon(t,x,U^\varepsilon)
=
0,
\qquad
U^\varepsilon(0,x,\theta)=U_0(x,\theta),
\end{equation}
where $y=(x,\theta)$ and
\begin{equation}\label{eq:exact-profile-flux}
G^\varepsilon(t,x,V)
:=
\left(
\frac{
F(\underline u+\varepsilon V)-F(\underline u)
}{\varepsilon}
\,;\,
\varepsilon^{-p}
\int_0^V
\bigl(
a(\underline u+\varepsilon s)-a(\underline u)
\bigr)
\cdot\nabla_x\phi\,ds
\right).
\end{equation}

\medskip
\noindent
On the other hand, the truncated profile
$U_{[K]}^\varepsilon$ defined in
\eqref{eq:truncated-amplitude-profile} satisfies the same
$\varepsilon$-dependent profile equation up to the residual produced
by the truncation:
\begin{equation}\label{eq:truncated-profile-equation}
\partial_tU_{[K]}^\varepsilon
+
\operatorname{div}_y
G^\varepsilon(t,x,U_{[K]}^\varepsilon)
=
\varepsilon^K R_K^\varepsilon,
\end{equation}
where $R_K^\varepsilon$ is uniformly bounded on the considered time
interval. 
Thus $U^\varepsilon$ and $U_{[K]}^\varepsilon$ have the same initial
datum and satisfy the same $\varepsilon$-dependent profile equation,
the only difference being a residual of order $\varepsilon^K$ in the
equation for $U_{[K]}^\varepsilon$. It remains to compare these two
profiles.

\paragraph{Stability and remainder estimate.}
\label{par:profile-stability}

\noindent
The exact profile $U^\varepsilon$ and the truncated profile
$U_{[K]}^\varepsilon$ have the same initial datum and satisfy the same
$\varepsilon$-dependent profile equation up to a residual of order
$\varepsilon^K$, uniformly with the derivatives required by the local
smooth theory. Moreover, the compatibility conditions imply that the
coefficients of the profile equation and their corresponding derivatives
remain uniformly bounded with respect to $\varepsilon$. Therefore,
standard local stability estimates for smooth solutions of quasilinear
scalar conservation laws \cite{PontryaginODE} yield, after possibly reducing the common time
interval,
\begin{equation}\label{eq:profile-W1infty-stability}
\left\|
U^\varepsilon-U_{[K]}^\varepsilon
\right\|_{W^{1,\infty}}
\lesssim
\varepsilon^K.
\end{equation}

\medskip
\noindent
We then define
\begin{equation}\label{eq:physical-remainder-definition}
r^\varepsilon(t,x)
:=
\varepsilon^{-K}
\left(
U^\varepsilon-U_{[K]}^\varepsilon
\right)
\left(
t,x,\frac{\phi(t,x)}{\varepsilon^p}
\right).
\end{equation}
Using \eqref{eq:physical-profile-difference}, we obtain
\[
u^\varepsilon(t,x)
=
\underline u(t,x)
+
\sum_{j=1}^K
\varepsilon^j
U^j\!\left(
t,x,\frac{\phi(t,x)}{\varepsilon^p}
\right)
+
\varepsilon^{K+1}r^\varepsilon(t,x).
\]
Estimate \eqref{eq:profile-W1infty-stability} immediately gives
\[
\|r^\varepsilon\|_{L^\infty(\Omega)}
\lesssim1.
\]
Furthermore, the chain rule gives, for
$z=t,x_1,\ldots,x_d$,
\[
\varepsilon^p\partial_z r^\varepsilon
=
\varepsilon^{p-K}
\partial_z
\left(U^\varepsilon-U_{[K]}^\varepsilon\right)
+
\varepsilon^{-K}
(\partial_z\phi)
\partial_\theta
\left(U^\varepsilon-U_{[K]}^\varepsilon\right),
\]
evaluated at $\theta=\phi/\varepsilon^p$. Hence, again by
\eqref{eq:profile-W1infty-stability},
\[
\|\varepsilon^p\nabla_{t,x}r^\varepsilon\|_{L^\infty(\Omega)}
\lesssim1,
\]
which proves \eqref{eq:complete-WKB-remainder}.

\subsection{The necessary condition}
\label{subsec:necessary-condition}

\noindent
We now prove the converse implication in
Theorem~\ref{thm:complete-WKB-criterion}. Assume that a complete
non-trivial WKB expansion of the form
\eqref{eq:complete-WKB-expansion}, with the remainder estimates
\eqref{eq:complete-WKB-remainder}, exists on $\Omega$. The argument is
based on the multidimensional Oleinik-type estimate
\begin{equation}\label{eq:oleinik-necessary}
\operatorname{div}_x a\bigl(u^\varepsilon(t,\cdot)\bigr)
\leq \frac{1}{t},
\qquad t>0,
\end{equation}
see, for instance, Hoff~\cite{Hoff83}. In the present smooth
setting, this estimate follows directly from the characteristic
equations and does not require the isotropy assumption arising in the
entropy-solution framework. Since the right-hand side of
\eqref{eq:oleinik-necessary} is independent of $\varepsilon$, the
singular terms generated by the fast variable in the WKB expansion
cannot persist with an arbitrary sign. We shall examine these terms
successively and show, using the periodicity and non-triviality of the
principal profile, that the compatibility conditions
\eqref{eq:propagated-generators} are necessarily satisfied. Once these
relations are recovered, the passage to the compatibility algebra
excludes the full-rank case $N_p=d$.

\subsubsection{Recovery of the compatibility conditions}
\label{subsubsec:necessary-recovery}

\noindent
Fix a time $t>0$ and choose once and for all $K\geq p$ in
\eqref{eq:complete-WKB-expansion}. The weighted derivative estimate in
\eqref{eq:complete-WKB-remainder} ensures that the remainder is
negligible at all the singular scales considered below. We first
examine explicitly the most singular contribution, which yields the
condition associated with $L_1$, and then repeat the same argument
inductively for $L_2,\ldots,L_{p-1}$.

\paragraph{The first compatibility condition.}
\label{par:necessary-first-compatibility}

\noindent
Write
\[
u^\varepsilon
=
\underline u+V^\varepsilon
+\varepsilon^{K+1} \, r^\varepsilon,
\qquad
V^\varepsilon
:=
\sum_{j=1}^{K}
\varepsilon^j \, 
U^j\!\left(
t,x,\frac{\phi(t,x)}{\varepsilon^p}
\right).
\]
Since
\[
\operatorname{div}_x a(u^\varepsilon)
=
a'(u^\varepsilon)\cdot\nabla_xu^\varepsilon,
\]
the singular terms come from differentiation of the fast variable
$\phi/\varepsilon^p$. In particular,
\[
\nabla_xV^\varepsilon
=
\sum_{j=1}^{K}
\varepsilon^j \, \nabla_x U^j
+
\sum_{j=1}^{K}
\varepsilon^{j-p} \, 
\nabla_x\phi \, \partial_\theta U^j,
\]
where the profiles are evaluated at
$\theta=\phi/\varepsilon^p$. From
\eqref{eq:complete-WKB-expansion} and
\eqref{eq:complete-WKB-remainder}, we have
\[
u^\varepsilon-\underline u
=
\varepsilon U_1+O(\varepsilon^2)
=
O(\varepsilon)
\]
uniformly on $\Omega$. Since $a$ is smooth, Taylor's formula gives
\[
a'(u^\varepsilon)
=
a'(\underline u)+O(\varepsilon)
\]
uniformly on $\Omega$. Consequently, the
most singular contribution to
$\operatorname{div}_x a(u^\varepsilon)$ is
\[
\varepsilon^{1-p} \, 
\bigl(
a'(\underline u)\cdot\nabla_x\phi
\bigr) \, 
\partial_\theta U^1
=
\varepsilon^{1-p} \, 
\bigl(L_1\phi\bigr) \, \partial_\theta U^1.
\]
All other contributions are of strictly lower singularity at this
stage. In particular, by
\eqref{eq:complete-WKB-remainder},
\[
\varepsilon^{K+1}\nabla_xr^\varepsilon
=
O(\varepsilon^{K+1-p}),
\]
which becomes negligible after multiplication by
$\varepsilon^{p-1}$ since $K\geq p$.

\medskip
\noindent
Multiplying \eqref{eq:oleinik-necessary} by
$\varepsilon^{p-1}$ and letting $\varepsilon\to0$, we therefore obtain
\begin{equation}\label{eq:necessary-first-sign}
\bigl(L_1\phi\bigr)(t,x)\,
\partial_\theta U^1(t,x,\theta)
\leq0.
\end{equation}
To interpret this inequality for an arbitrary value of $\theta$, we use a sampling argument. If $\phi(t,x)\neq0$, then, for
every $\theta\in\mathbb T$, one can choose a sequence
$\varepsilon_\nu\to0$ such that
\[
\frac{\phi(t,x)}{\varepsilon_\nu^p}
\longrightarrow\theta
\qquad\text{in }\mathbb T.
\]
Passing to this sequence in the Oleinik inequality shows that
\eqref{eq:necessary-first-sign} holds for every
$\theta\in\mathbb T$. If $\phi(t,x)=0$, the non-triviality condition
\eqref{eq:nontrivial-phase-p} gives
$\nabla_x\phi(t,x)\neq0$; hence $(t,x)$ can be approached by points
$(t,x_\nu)$ with $\phi(t,x_\nu)\neq0$, and the same conclusion follows
by continuity. Thus, for every $(t,x)\in\Omega$ and every $\theta\in\mathbb T$,
\[
\bigl(L_1\phi\bigr)(t,x)
\,\partial_\theta U^1(t,x,\theta)
\leq0.
\]
Since $L_1\phi$ is independent of $\theta$, if
$(L_1\phi)(t,x)\neq0$, then
$\partial_\theta U^1(t,x,\cdot)$ has a constant sign on
$\mathbb T$. Periodicity forces it to vanish identically, in
contradiction with the non-triviality condition
\eqref{eq:nontrivial-profile-p}. Therefore
\begin{equation}\label{eq:necessary-L1}
L_1(t)\phi(t,\cdot)=0
\qquad\text{on }\Omega.
\end{equation}

\paragraph{Induction on the compatibility conditions.}
\label{par:necessary-induction}

\noindent
Assume now that, for some $q\in\{2,\ldots,p-1\}$,
\begin{equation}\label{eq:necessary-induction-hypothesis}
L_k(t)\phi(t,\cdot)=0,
\qquad
k=1,\ldots,q-1.
\end{equation}
Expanding $a'(u^\varepsilon)$ around $\underline u$, the induction
hypothesis cancels all the preceding singular contributions. The first
remaining one is
\[
\frac{\varepsilon^{q-p}}{q!} \, 
\bigl(
a^{(q)}(\underline u)\cdot\nabla_x\phi
\bigr) \, \partial_\theta\bigl( (U^1)^q\bigr) 
=
\frac{\varepsilon^{q-p}}{q!} \, 
\bigl(L_q\phi\bigr) \, \partial_\theta\bigl( (U^1)^q\bigr) .
\]
The remainder is again negligible at this scale, since after
multiplication by $\varepsilon^{p-q}$ its contribution is
$O(\varepsilon^{K+1-q})$. Multiplying
\eqref{eq:oleinik-necessary} by $\varepsilon^{p-q}$, passing to the
limit and using the same sampling argument as above, we obtain
\begin{equation}\label{eq:necessary-q-sign}
\bigl(L_q\phi\bigr)(t,x) \, \partial_\theta\bigl( (U^1)^q\bigr) (t,x,\theta)
\leq0,
\qquad
\forall \, \theta\in\mathbb T.
\end{equation}
If $(L_q\phi)(t,x)\neq0$, periodicity implies that
$\partial_\theta\bigl( (U^1)^q\bigr) (t,x,\cdot)$ vanishes identically. Hence
$ (U^1)^q(t,x,\cdot)$ is constant. Since the function $U^1$ is real-valued and
continuous, $U^1(t,x,\cdot)$ must be constant, contradicting
\eqref{eq:nontrivial-profile-p}. Thus
\[
L_q(t)\phi(t,\cdot)=0.
\]
Starting from \eqref{eq:necessary-L1} and iterating this argument for
$q=2,\ldots,p-1$, we conclude that
\begin{equation}\label{eq:necessary-all-compatibilities}
L_k(t)\phi(t,\cdot)=0,
\qquad
k=1,\ldots,p-1.
\end{equation}

\subsubsection{From the generators to the  Lie algebra}
\label{subsubsec:necessary-lie-algebra}

\noindent
The preceding argument has recovered all the compatibility conditions
associated with the generators $L_1,\ldots,L_{p-1}$. It remains to
translate these relations into a statement on the whole Lie 
algebra. As in the proof of Lemma~\ref{Lemma propagation Lie algebra}, the
property of annihilating $\phi(t,\cdot)$ is preserved under smooth
linear combinations and Lie brackets. Therefore
\[
T\phi(t,\cdot)=0,
\qquad
\forall \, T\in
\operatorname{Lie}\bigl(L_1(t),\ldots,L_{p-1}(t)\bigr).
\]
By the definition of $\mathcal A_p$, this is equivalent to
\begin{equation}\label{eq:necessary-Ap-kernel}
\mathcal A_p(t,x)
\subset
\ker d_x\phi(t,x),
\qquad
(t,x)\in\Omega.
\end{equation}

\subsubsection{Exclusion of the full-rank case}
\label{subsubsec:necessary-full-rank}

\noindent
We can conclude using the geometric condition obtained above.
Suppose, by contradiction, that
$ N_p=d $. By Lemma~\ref{lem:propagation-rank-Ap}, the compatibility algebra has
full rank throughout $\Omega$, and therefore
\[
\mathcal A_p(t,x)=T_x\mathbb R^d,
\qquad
(t,x)\in\Omega.
\]
Together with \eqref{eq:necessary-Ap-kernel}, this gives
\[
d_x\phi(t,x)=0,
\qquad
(t,x)\in\Omega,
\]
or equivalently $ \nabla_x\phi(t,x)=0 $.
This contradicts the non-triviality condition
\eqref{eq:nontrivial-phase-p}. Hence $
N_p<d $. This proves the necessary implication and completes the proof of
Theorem~\ref{thm:complete-WKB-criterion}.

\subsection{Explicit geometric criteria for compatible phases}
\label{subsec:explicit-geometric-criteria}

\noindent
The necessary condition obtained above shows that any non-trivial phase
associated with a complete supercritical WKB expansion must satisfy
\[
L_j\phi=0,
\qquad
j=1,\ldots,p-1.
\]
Equivalently, the phase must be constant along the directions generated
by the compatibility vector fields. The abstract criterion in
Theorem~\ref{thm:complete-WKB-criterion} is therefore expressed in terms
of the Lie algebra generated by these fields and its rank $N_p$. In several relevant situations, this geometric condition can be made
explicit. The bracket structure of the vector fields allows one to
determine when new independent directions are generated and, conversely,
when the compatibility distribution closes. We record below some
representative cases.

\subsubsection{Lie-bracket structure of the compatibility fields}
\label{subsubsec:Lie-bracket-structure}

\noindent
We first record a general identity for the Lie brackets of the
compatibility vector fields. It will be the basic tool in the explicit
geometric analysis carried out below.

\begin{lem}\label{Lie bracket}
For every $i,j\in\{1,\ldots,p-1\}$,
\begin{equation}\label{eq:bracket-identity}
[L_i,L_j]
=
(L_i\underline u_0)L_{j+1}
-
(L_j\underline u_0)L_{i+1}.
\end{equation}
\end{lem}

\begin{Proof}
Let $v$ be a smooth function. By definition,
\[
L_i
=
a^{(i)}(\underline u_0)\cdot\nabla_x,
\qquad
L_j
=
a^{(j)}(\underline u_0)\cdot\nabla_x.
\]
Using the chain rule, we obtain
\begin{align*}
L_iL_jv
&=
a^{(i)}(\underline u_0)\cdot
\nabla_x
\left(
a^{(j)}(\underline u_0)\cdot\nabla_xv
\right)
\\
&=
\bigl(L_i\underline u_0\bigr)
a^{(j+1)}(\underline u_0)\cdot\nabla_xv
+
D_x^2v
\bigl(
a^{(i)}(\underline u_0),
a^{(j)}(\underline u_0)
\bigr)
\\
&=
\bigl(L_i\underline u_0\bigr)L_{j+1}v
+
D_x^2v
\bigl(
a^{(i)}(\underline u_0),
a^{(j)}(\underline u_0)
\bigr).
\end{align*}
Similarly,
\[
L_jL_iv
=
\bigl(L_j\underline u_0\bigr)L_{i+1}v
+
D_x^2v
\bigl(
a^{(j)}(\underline u_0),
a^{(i)}(\underline u_0)
\bigr).
\]
Since the Hessian $D_x^2v$ is symmetric, the second-order terms cancel.
Hence,
\[
[L_i,L_j]v
=
\bigl(L_i\underline u_0\bigr)L_{j+1}v
-
\bigl(L_j\underline u_0\bigr)L_{i+1}v.
\]
Since this identity holds for every smooth function $v$, we obtain
\cref{eq:bracket-identity}.
\end{Proof}

\subsubsection{The linearly nondegenerate case}
\label{par:linearly-nondegenerate-case}

\noindent
We first consider the case in which the successive derivatives of the
characteristic velocity remain linearly independent. In this situation,
the higher-order vector fields generated by the bracket identity
\cref{eq:bracket-identity} cannot, in general, be absorbed into the
compatibility distribution $\mathcal D_p$ defined in
\cref{eq:distribution-Dp}. Thus the Lie closure $\mathcal A_p$ may be
strictly larger than $\mathcal D_p$. The involutivity condition therefore
imposes additional constraints on the background state $\underline u_0$.

\medskip
\noindent
To formulate this case precisely, we use the following notion of linear
nondegeneracy.

\begin{defn}[Linear nondegeneracy up to order $j$]
\label{linear nondegeneracy}
Let $1\leq j\leq d$. A flux
$F\in C^\infty(\mathbb{R};\mathbb{R}^d)$ is said to be
\emph{linearly nondegenerate up to order $j$} if
\[
\operatorname{rank}
\bigl\{
a'(u),a''(u),\ldots,a^{(j)}(u)
\bigr\}
=j,
\qquad
\forall u\in\mathbb{R}.
\]
\end{defn}

\smallskip
\noindent
A simple example which is linearly nondegenerate up to order $d$ is
given by the velocity field
\begin{equation}\label{example linear nondegenerate flux}
a(u)=\left(
\frac{u}{1!},
\frac{u^2}{2!},
\ldots,
\frac{u^d}{d!}
\right).
\end{equation}

\begin{prop}
\label{prop:compatible-phase-p3}
Assume that $d=p=3$ and that the velocity field $a$ is given by
\cref{example linear nondegenerate flux}. Then there exists locally a
nonconstant phase $\phi_0$ satisfying 
$L_1\phi_0=0$ and $L_2\phi_0=0$
if and only if
\begin{equation}\label{eq:condition-compatible-phase-p3}
L_1\underline u_0=0.
\end{equation}
\end{prop}

\begin{Proof}
For $p=3$, the compatibility distribution \cref{eq:distribution-Dp} is
\[
\mathcal D
=
\operatorname{span}\{L_1,L_2\}.
\]
For the velocity field \cref{example linear nondegenerate flux}, one has
\[
\det\bigl(
a'(\underline u_0),
a''(\underline u_0),
a'''(\underline u_0)
\bigr)
=
1.
\]
Hence $L_1,L_2,L_3$ are pointwise linearly independent, and in particular
$\mathcal D$ has constant rank two.
By the Frobenius theorem, a local nonconstant phase satisfying the
compatibility system exists if and only if $\mathcal D$ is involutive.
Since $\mathcal D$ is generated by $L_1$ and $L_2$, this is equivalent to
$[L_1,L_2]\in\mathcal D$.
Using \cref{eq:bracket-identity}, we have
\[
[L_1,L_2]
=
(L_1\underline u_0)L_3
-
(L_2\underline u_0)L_2.
\]
Moreover,
\begin{equation}\label{eq:det-generic-rank}
\det\bigl(L_1,L_2,[L_1,L_2]\bigr)
=
(L_1\underline u_0)
\det\bigl(L_1,L_2,L_3\bigr).
\end{equation}
Since $L_1,L_2,L_3$ are pointwise linearly independent,
\cref{eq:det-generic-rank} shows that
$[L_1,L_2]\in\mathcal D$ if and only if
$L_1\underline u_0=0$.
Under \cref{eq:condition-compatible-phase-p3}, the bracket identity reduces
to
\[
[L_1,L_2]
=
-(L_2\underline u_0)L_2
\in\mathcal D,
\]
so $\mathcal D$ is involutive and therefore
$\mathcal A=\mathcal D$. Since this distribution has rank two, the
Frobenius theorem yields locally a nonconstant phase satisfying the
compatibility system, with one transverse degree of freedom.
\end{Proof}

\medskip
\noindent
The previous argument extends to arbitrary $p\leq d$ when the vectors
$a'(u),\ldots,a^{(p)}(u)$ remain linearly independent. In that case, the
only possible obstruction to the involutivity of the compatibility
distribution comes from the $L_p$ component appearing in the brackets
with $L_{p-1}$.

\begin{prop}
\label{prop:general-compatible-phase}
Assume that $p\leq d$ and that $
a'(u),a''(u),\ldots,a^{(p)}(u) $ 
are linearly independent for every $u\in\mathbb{R}$. Then, the distribution $\mathcal D$ defined in
\cref{eq:distribution-Dp} is involutive if and only if
\begin{equation}\label{eq:general-involutivity-condition}
L_j\underline u_0=0,
\qquad
j=1,\ldots,p-2.
\end{equation}
\end{prop}

\begin{Proof}
The linear independence of $
a'(\underline u_0),\ldots,a^{(p)}(\underline u_0) $ implies that the vector fields $ L_1,\ldots,L_p $ are pointwise linearly independent. In particular $
\operatorname{rank}\mathcal D=p-1 $ and 
$ L_p\notin\mathcal D $. We now use the bracket identity \cref{eq:bracket-identity}. If $ 1\leq i,j\leq p-2 $, then $ 
L_{i+1},L_{j+1}\in\mathcal D $, 
and therefore $ [L_i,L_j]\in\mathcal D $.
It remains, up to antisymmetry of the Lie bracket, to consider
\[
[L_i,L_{p-1}],
\qquad
i=1,\ldots,p-2.
\]
By \cref{eq:bracket-identity},
\[
[L_i,L_{p-1}]
=
(L_i\underline u_0)L_p
-
(L_{p-1}\underline u_0)L_{i+1}.
\]
The second term belongs to $\mathcal D$, whereas
$L_p\notin\mathcal D$. Consequently,
\[
[L_i,L_{p-1}]\in\mathcal D
\quad\Longleftrightarrow\quad
L_i\underline u_0=0.
\]
Hence \cref{eq:general-involutivity-condition} holds if and only if
$\mathcal D$ is involutive, in which case $\mathcal A=\mathcal D$ and
$\operatorname{rank}\mathcal D=p-1$. By the Frobenius theorem, it is
locally tangent to leaves of dimension $p-1$, leaving $d-p+1$
transverse degrees of freedom for the construction of compatible phases.
\end{Proof}

\subsubsection{Closure of the derivative hierarchy}
\label{par:closure-derivative-hierarchy}

\noindent
We now consider the complementary situation in which the successive
derivatives of the characteristic velocity cease to generate new
independent directions. More precisely, we assume that the derivative
hierarchy closes at some finite order $j$. In this case, all higher-order
vector fields remain in the span of the first $j-1$ fields, and the bracket
identity \cref{eq:bracket-identity} implies involutivity without imposing
additional conditions on the background state $\underline u_0$.

\begin{prop}
\label{prop:closure-at-order-j}
Assume that, for some $j\geq2$, the directions
$a'(u),a''(u),\ldots,a^{(j-1)}(u)$
are linearly independent for every $u\in\mathbb{R}$, and that
\[
a^{(j)}(u)
\in
\operatorname{span}
\bigl\{
a'(u),a''(u),\ldots,a^{(j-1)}(u)
\bigr\},
\qquad
\forall u\in\mathbb{R}.
\]
Then, for every $p\geq j$, the compatibility distribution
$\mathcal D=\operatorname{span}\{L_1,\ldots,L_{j-1}\}$ is involutive,
and consequently $\mathcal A=\mathcal D$.
\end{prop}

\begin{Proof}
By the linear independence assumption, there exist unique smooth functions
$c_1,\ldots,c_{j-1}$ such that
\[
a^{(j)}(u)
=
\sum_{m=1}^{j-1}
c_m(u)\,a^{(m)}(u).
\]
Differentiating with respect to $u$ gives
\[
a^{(j+1)}(u)
=
\sum_{m=1}^{j-1}
c_m'(u)\,a^{(m)}(u)
+
\sum_{m=1}^{j-1}
c_m(u)\,a^{(m+1)}(u).
\]
Every term on the right-hand side belongs to
\[
\operatorname{span}
\bigl\{
a'(u),\ldots,a^{(j-1)}(u)
\bigr\},
\]
because the only derivative of order $j$ which occurs is $a^{(j)}$,
which belongs to this span by assumption. Hence
\[
a^{(j+1)}(u)
\in
\operatorname{span}
\bigl\{
a'(u),\ldots,a^{(j-1)}(u)
\bigr\}.
\]
Since the coefficients in such a decomposition are again smooth, the same
argument can be iterated. Thus, for every $k\geq j$,
\[
a^{(k)}(u)
\in
\operatorname{span}
\bigl\{
a'(u),\ldots,a^{(j-1)}(u)
\bigr\}.
\]

\medskip
\noindent
Evaluating at $u=\underline u_0(x)$ yields, for every $k\geq j$,
\[
L_k
\in
\operatorname{span}\{L_1,\ldots,L_{j-1}\}.
\]
Consequently, for every $p\geq j$,
\[
\mathcal D
=
\operatorname{span}\{L_1,\ldots,L_{p-1}\}
=
\operatorname{span}\{L_1,\ldots,L_{j-1}\}.
\]

\medskip
\noindent
To prove involutivity, let $1\leq r,s\leq j-1$. By
\cref{eq:bracket-identity},
\[
[L_r,L_s]
=
(L_r\underline u_0)L_{s+1}
-
(L_s\underline u_0)L_{r+1}.
\]
Since $r+1,s+1\leq j$ and
\[
L_j\in\operatorname{span}\{L_1,\ldots,L_{j-1}\},
\]
both terms on the right-hand side belong to
$\operatorname{span}\{L_1,\ldots,L_{j-1}\}$. Hence $\mathcal D$ is
involutive, and therefore $\mathcal A=\mathcal D$.

\medskip
\noindent
Finally, the linear independence of
$a',\ldots,a^{(j-1)}$ implies that
$\operatorname{rank}\mathcal D=\operatorname{rank}\mathcal A=j-1$.
By the Frobenius theorem, the corresponding leaves have dimension $j-1$,
leaving $d-j+1$ transverse degrees of freedom for compatible phases.
\end{Proof}

\medskip
\noindent
These examples illustrate how the abstract condition $N_p<d$ can be
related to explicit structural properties of the flux and the background
state. In the situations considered here, the compatibility distribution
$\mathcal D_p$ is involutive, and therefore coincides with its Lie closure
$\mathcal A_p$.

\subsection{Cancellation under first-order incompatibility}
\label{Cancellation-p>2}

\noindent
We now prove Theorem~\ref{thm-cancel3}. The argument follows the same
mechanism as in the cancellation regime of Section~\ref{sec:2}. The
higher oscillatory frequency $\varepsilon^{-p}$ changes the natural
boundary-layer scale: the appropriate time-space localization is now of
order $\varepsilon^{p-1}$. After rescaling at this scale, the localized
problem is reduced to the canonical cancellation regime already treated
for $p=2$.

\subsubsection{Rescaling and reduction to the canonical cancellation regime}
\label{cancellation-rescaling-p>2}

\noindent
We first identify the scale on which the supercritical oscillations
enter the canonical cancellation regime. Let $k\in\mathbb Z^d$ index
the hypercubes of the spatial localization and introduce the rescaled
variables
\begin{equation}\label{eq:cancellation-rescaling-p>2}
t=\varepsilon^{p-1}\tau,
\qquad
x=\varepsilon^{p-1}(k+y).
\end{equation}
We define
\begin{equation}\label{eq:rescaled-solution-cancellation-p>2}
v_k^\varepsilon(\tau,y)
:=
u^\varepsilon\!\left(
\varepsilon^{p-1}\tau,
\varepsilon^{p-1}(k+y)
\right).
\end{equation}
Since
\[
\partial_tu^\varepsilon
=
\varepsilon^{-(p-1)}
\partial_\tau v_k^\varepsilon,
\qquad
\nabla_x
=
\varepsilon^{-(p-1)}\nabla_y,
\]
the scalar conservation law \cref{eq:claw} keeps the same form:
\begin{equation}\label{eq:rescaled-claw-cancellation-p>2}
\partial_\tau v_k^\varepsilon
+
\operatorname{div}_yF(v_k^\varepsilon)
=
0.
\end{equation}
At $\tau=0$, the corresponding initial data are
\[
v_k^\varepsilon(0,y)
=
\underline u_0\!\left(
\varepsilon^{p-1}(k+y)
\right)
+
\varepsilon
U_0\!\left(
\varepsilon^{p-1}(k+y),
\frac{
\phi_0\!\left(
\varepsilon^{p-1}(k+y)
\right)
}{\varepsilon^p}
\right).
\]

\medskip
\noindent
Expanding the phase around the center
$\varepsilon^{p-1}k$ of the hypercube gives
\[
\phi_0\!\left(
\varepsilon^{p-1}(k+y)
\right)
=
\phi_0(\varepsilon^{p-1}k)
+
\varepsilon^{p-1}
\nabla_x\phi_0(\varepsilon^{p-1}k)\cdot y
+
O(\varepsilon^{2p-2}),
\]
and therefore
\begin{equation}\label{eq:rescaled-phase-cancellation-p>2}
\frac{
\phi_0\!\left(
\varepsilon^{p-1}(k+y)
\right)
}{\varepsilon^p}
=
\frac{
\phi_0(\varepsilon^{p-1}k)
}{\varepsilon^p}
+
\frac{
\nabla_x\phi_0(\varepsilon^{p-1}k)\cdot y
}{\varepsilon}
+
O(\varepsilon^{p-2}).
\end{equation}
Thus, in the rescaled variables, the oscillations have frequency
$\varepsilon^{-1}$. Moreover, Assumption~\ref{ass-fullyn}
ensures that the interaction between the phase and the first
characteristic field remains uniformly non-degenerate. The localized
problem is therefore reduced to the canonical cancellation regime of
Section~\ref{sec:2}.

\subsubsection{Local cancellation at the boundary-layer time}
\label{cancellation-local-p>2}

\noindent
We may now apply the cancellation mechanism developed in
Section~\ref{sec:2}. In the rescaled variables, the weakly nonlinear
geometric-optics reduction leads to the same Burgers-type dynamics in
the fast variable, and the corresponding Oleinik estimate yields decay
of the oscillatory component.

\medskip
\noindent
More precisely, the analogue of
Proposition~\ref{prop:local-cancellation} gives, at
$\tau=\mathcal I$,
\begin{equation}\label{eq:local-cancellation-p>2}
\left\|
v_k^\varepsilon(\mathcal I,\cdot)
-
\underline v_k^\varepsilon(\mathcal I,\cdot)
-
\varepsilon\overline U_k^\varepsilon(\mathcal I,\cdot)
\right\|_{L^1(B(0,1))}
\leq
C_\rho\frac{\varepsilon}{\mathcal I}.
\end{equation}
Here $\overline U_k^\varepsilon$ denotes the averaged profile associated
with the localized problem. As in Section~\ref{sec:2}, the constants
are uniform with respect to the hypercube index $k$ as long as the
corresponding hypercubes remain in the localization region.

\subsubsection{Return to the physical variables}
\label{cancellation-return-p>2}

\noindent
It remains to transfer the local estimate
\eqref{eq:local-cancellation-p>2} to the physical variables and to
propagate it to later times. Compared with the case $p=2$, the only
change is the scale of localization: the physical hypercubes now have
size $\varepsilon^{p-1}$, and the end of the initial boundary layer
corresponds to $ t=\varepsilon^{p-1}\mathcal I $. The same reconstruction and covering argument as in
Section~\ref{sec:2} therefore yields the corresponding global
$L^1$ estimate at this intermediate time.

\medskip
\noindent
Finally, the local $L^1$-contraction argument of
Section~\ref{sec:2} propagates this estimate to later times. Choosing
first $\mathcal I$ sufficiently large and then $\varepsilon$
sufficiently small, we obtain
\[
u^\varepsilon
=
\underline u
+
\varepsilon\overline U
+
o(\varepsilon)
\]
in the sense of the $L^1$-norm. This is precisely
\eqref{eq:expansion-avg-p>2} and completes the proof of
Theorem~\ref{thm-cancel3}.


\appendix
 
\section{Regularity and compactness of the nonlinear semigroup}
\label{app:regularity-compactness-semigroup}

\medskip
\noindent The propagation or lack of propagation of high frequencies is closely related to the optimal regularity of entropy solutions. More precisely, it is linked to the compactness properties of the Kružkov operator, that is, to uniform bounds on the regularity gained by entropy solutions. More generally, the regularity theory for entropy solutions is a broad topic, For approaches based on the kinetic formulation,  see \cite{COW,JabinX,LPT}. Other several approaches have been developed \cite{ACCM26,Marconi,Sil}, and small-scale effects have also been studied in different contexts \cite{Dal09,EE}. In this section, we shall focus mainly on the maximal compactness properties of the Kružkov operator. 
We investigate here the fractional Sobolev regularity of the highly
oscillating solutions constructed in the propagation regime of
Section~\ref{sec:3}, and relate their WKB description to the regularity and
compactness properties of the nonlinear semigroup associated with
\cref{eq:claw}. Under the assumptions of Theorem~\ref{thm:complete-WKB-criterion}, the exact
solution admits, for every fixed $K\geq1$, the expansion
\[
u^\varepsilon=v_K^\varepsilon+\varepsilon^{K+1}r^\varepsilon,
\]
where $v_K^\varepsilon$ is the truncated WKB expansion defined in
\eqref{eq:truncated-WKB}, and the remainder satisfies
\cref{eq:complete-WKB-remainder}. Since the
principal oscillation has amplitude $\varepsilon$ and wavelength of order
$\varepsilon^p$, the expected fractional Sobolev scaling is
$\varepsilon(\varepsilon^p)^{-s}=\varepsilon^{1-sp}$, with critical exponent
$s_c=1/p$. The following proposition makes this precise.

\begin{prop}[Critical fractional Sobolev scaling]
\label{prop:critical-sobolev-scaling}
Assume the hypotheses of
Theorem~\ref{thm:complete-WKB-criterion}, and let $K\geq1$ be fixed.
For every fixed $t\in[0,T]$, there exists a ball
$B\Subset\Omega_t$ such that, for every $0<s<1$, there exist constants
$c_s,C_s>0$ and $\varepsilon_s>0$ satisfying 
(for every $0<\varepsilon\leq\varepsilon_s$)
\begin{equation}\label{eq:critical-sobolev-scaling}
c_s\,\varepsilon^{1-sp}
\leq
\|v_K^\varepsilon(t,\cdot)-\underline u(t,\cdot)\|_{W^{s,1}(B)}
\leq
C_s\,\varepsilon^{1-sp} .
\end{equation}
\end{prop}

\noindent 
Thus $s_c =1/p$ is the critical fractional Sobolev exponent: the oscillatory
perturbation converges strongly to zero in $W^{s,1}(B)$ for $s<1/p$, remains
bounded above and away from zero for $s=1/p$, and its $W^{s,1}(B)$ norm
diverges for $s>1/p$. This means that the semi-group cannot induce a regularizing $ W^{s,1} $-effect better than $ 1/p $.

\smallskip

\noindent The question is: {\it what is the largest $p$ for which high-frequency waves can be propagated$\, $?}

\smallskip

\noindent We adopt the convention $ N_1 = 0 $. Our study shows that the largest $p$ is given by\footnote{Remark that calculating $ p_{\max} $ relies on the construction of the Lie closure $ \mathcal A_p $ for $ p \leq p_{\max} $.}
\[
p_{\max}:=\max\{p\geq 1:\,N_p<d\} . 
\]
The larger $p_{\max}$ is, the weaker the  uniform regularizing effect is. 
For a linear flux $ F $, we find $ p_{\max} = + \infty $ and, as expected for a linear equation, all Sobolev regularities are simply propagated. As soon as  $ a' \not \equiv 0 $, we can find $ u_0 $ such that $ a' \circ u_0 \not \equiv 0 $ so that $ N_2 = 1 $ is reachable. For $ d= 1 $, we recover $ p_{\max} = 1 $ which corresponds to $ BV $. For $ 1 < d $ and a nonlinear flux, we obtain $ 2 \leq p_{\max} < + \infty $ (see \cite{J4} for the right hand side inequality) so that the fractional regularity $ s $ of the semi-group must be adjusted in such a way that 
$$ s \leq s_{max} := 1 / p_{max} \leq 1/2 \quad \text{when} \quad 1 < d .$$
In \eqref{eq:critical-sobolev-scaling}, the upper estimate follows from the same amplitude--frequency argument as in
\cite{J4}: the perturbation has size $\varepsilon$ in $L^1$, while one
spatial derivative of the fast variable costs a factor $\varepsilon^{-p}$,
and interpolation  \cite{T2} gives $\varepsilon^{1-sp}$. For the lower estimate, we use the translation representation (Gagliardo seminorm)
\begin{equation}\label{eq:appendix-gagliardo-translations}
[f]_{W^{s,1}(B)}
=
\int_{\mathbb R^d}
\frac{1}{|H|^{d+s}}
\int_{B\cap(B-H)}
|f(x+H)-f(x)|\,dx\,dH.
\end{equation}
Since the oscillations have wavelength of order $\varepsilon^p$, the natural
translations are $H=\varepsilon^p h$. The lower bound is based on the
following estimate.

\begin{lem}[Wavelength-scale increment]
\label{lem:appendix-wavelength-increment}
There exist balls $B_0\Subset B\Subset\Omega_t$, a compact set
$\mathcal H\Subset\mathbb R^d\setminus\{0\}$ with non-empty interior, and $c>0$ such that, for all sufficiently small $\varepsilon$,
\begin{equation}\label{eq:appendix-increment-lower}
\int_{B_0}
\left|
(v_K^\varepsilon-\underline u)(t,x+\varepsilon^p h)
-
(v_K^\varepsilon-\underline u)(t,x)
\right|
\,dx
\geq
c \, \varepsilon,
\qquad
h\in\mathcal H.
\end{equation}
\end{lem}

\medskip
\noindent
The appendix is organized in two parts. In
Subsection~\ref{subsec:appendix-proof-critical-scaling}, we prove
Proposition~\ref{prop:critical-sobolev-scaling} and transfer the resulting
scaling to the exact solutions and the entropy semigroup. In
Subsection~\ref{subsec:appendix-proof-wavelength-increment}, we prove
Lemma~\ref{lem:appendix-wavelength-increment}, first establishing the local
averaging property needed for the variable phase.


\subsection{Proof of Proposition~\ref{prop:critical-sobolev-scaling}
and consequences}
\label{subsec:appendix-proof-critical-scaling}
We first establish the two-sided fractional Sobolev estimate for the truncated
WKB expansion $v_K^\varepsilon$. The upper bound follows from the size of the
oscillation and its wavelength, while the lower bound relies on the
wavelength-scale increment estimate of
Lemma~\ref{lem:appendix-wavelength-increment}. We then transfer the same
scaling to the exact solution $u^\varepsilon$ and derive the corresponding
compactness consequences for the entropy semigroup.

\subsubsection{Proof of Proposition~\ref{prop:critical-sobolev-scaling}}
\label{subsubsec:appendix-proof-critical-scaling}

We begin with the upper bound. From \cref{eq:truncated-WKB} and the boundedness
of the profiles,
\begin{equation}\label{eq:appendix-upper-L1}
\|v_K^\varepsilon(t,\cdot)-\underline u(t,\cdot)\|_{L^1(B)}
\lesssim
\varepsilon.
\end{equation}
Moreover, for every $k=1,\ldots,K$,
\[
\begin{aligned}
\nabla_x
\left[
\varepsilon^k
U^k\!\left(
t,x,\frac{\phi(t,x)}{\varepsilon^p}
\right)
\right]
={}&
\varepsilon^k
\nabla_xU^k\!\left(
t,x,\frac{\phi(t,x)}{\varepsilon^p}
\right)
\\
&+
\varepsilon^{k-p}
\partial_\theta U^k\!\left(
t,x,\frac{\phi(t,x)}{\varepsilon^p}
\right)
\nabla_x\phi(t,x).
\end{aligned}
\]
The most singular contribution corresponds to $k=1$, hence
\begin{equation}\label{eq:appendix-upper-W11}
\|v_K^\varepsilon(t,\cdot)-\underline u(t,\cdot)\|_{W^{1,1}(B)}
\lesssim
\varepsilon^{1-p}.
\end{equation}
Interpolating between \cref{eq:appendix-upper-L1,eq:appendix-upper-W11} \cite{T2},
\[
\|v_K^\varepsilon(t,\cdot)-\underline u(t,\cdot)\|_{W^{s,1}(B)}
\lesssim
\varepsilon^{1-s}
\bigl(\varepsilon^{1-p}\bigr)^s
=
\varepsilon^{1-sp}.
\]
For the lower bound, restricting the translation variable in
\cref{eq:appendix-gagliardo-translations} to
$H\in\varepsilon^p\mathcal H$ and using
Lemma~\ref{lem:appendix-wavelength-increment}, we obtain
\[
[v_K^\varepsilon(t,\cdot)-\underline u(t,\cdot)]_{W^{s,1}(B)}
\gtrsim
\varepsilon
\int_{\varepsilon^p\mathcal H}
\frac{dH}{|H|^{d+s}} =
\varepsilon^{1-sp}
\int_{\mathcal H}
\frac{dh}{|h|^{d+s}}
\gtrsim
\varepsilon^{1-sp}.
\]
This proves \cref{eq:critical-sobolev-scaling}.

\subsubsection{Consequences for the exact solutions and the semigroup}
\label{subsubsec:appendix-exact-semigroup}

By \cref{eq:complete-WKB-expansion,eq:complete-WKB-remainder},
\[
\|u^\varepsilon(t,\cdot)-v_K^\varepsilon(t,\cdot)\|_{L^1(B)}
\lesssim
\varepsilon^{K+1},
\qquad
\|u^\varepsilon(t,\cdot)-v_K^\varepsilon(t,\cdot)\|_{W^{1,1}(B)}
\lesssim
\varepsilon^{K+1-p}.
\]
Interpolation gives
\begin{equation}\label{eq:appendix-exact-remainder}
\|u^\varepsilon(t,\cdot)-v_K^\varepsilon(t,\cdot)\|_{W^{s,1}(B)}
\lesssim
\varepsilon^{K+1-sp}
=
o\!\left(\varepsilon^{1-sp}\right).
\end{equation}
Hence, for every $0<s<1$, there exist $c_s',C_s'>0$ such that
\begin{equation}\label{eq:appendix-exact-scaling}
c_s'\,\varepsilon^{1-sp}
\leq
\|u^\varepsilon(t,\cdot)-\underline u(t,\cdot)\|_{W^{s,1}(B)}
\leq
C_s'\,\varepsilon^{1-sp},
\end{equation}
for all sufficiently small $\varepsilon$.

\medskip
\noindent
Consequently, for every fixed $t\in(0,T]$,
$u^\varepsilon(t,\cdot)\to\underline u(t,\cdot)$ strongly in $W^{s,1}(B)$
when $0<s<1/p$; at the critical exponent $s=1/p$, the norm
$\|u^\varepsilon(t,\cdot)-\underline u(t,\cdot)\|_{W^{1/p,1}(B)}$ remains
bounded above and away from zero; and for $s>1/p$ it tends to $+\infty$. Let $S_t$ denote the entropy semigroup associated with \cref{eq:claw}.
Since $S_tu_0^\varepsilon=u^\varepsilon(t,\cdot)$ and
$u^\varepsilon(t,\cdot)-\underline u(t,\cdot)\to0$ in $L^1(B)$, the lower
bound in \cref{eq:appendix-exact-scaling} at $s=1/p$ shows that the family
cannot be relatively compact in $W^{1/p,1}(B)$. Thus strong convergence
holds below the critical exponent, compactness is lost at the critical
exponent, and the fractional Sobolev norms become unbounded above it.

\subsection{Proof of the wavelength-scale increment estimate}
\label{subsec:appendix-proof-wavelength-increment}

We now prove Lemma~\ref{lem:appendix-wavelength-increment}, which provides
the key lower bound in Proposition~\ref{prop:critical-sobolev-scaling}.
The main difference with the wavelength-scale argument of \cite{J4} is that
the fast phase is variable. Indeed, Taylor's formula gives
\begin{equation}\label{eq:appendix-variable-phase-motivation}
\frac{
\phi(t,x+\varepsilon^p h)-\phi(t,x)
}{
\varepsilon^p
}
=
\nabla_x\phi(t,x)\cdot h
+
O(\varepsilon^p).
\end{equation}
A translation of size $\varepsilon^p h$ gives asymptotically
a translation in the variable $ \theta $. We first establish the corresponding
local averaging result, and then use it to prove the increment estimate.


\subsubsection{Local averaging for the variable phase}
\label{subsubsec:appendix-local-averaging}

\begin{lem}[Local averaging for a variable phase]
\label{lem:local-averaging-variable-phase}
Assume that $\phi$ is non-trivial in the sense of
\cref{eq:nontrivial-phase-p}, and fix $t\in[0,T]$.
For every $x_0\in\Omega_t$, there exists a ball
$B\Subset\Omega_t$ centered at $x_0$ such that, for every compact set
$\mathcal H\Subset\mathbb R^d$ and every function
$G\in C(\overline B\times\mathbb T\times\mathcal H)$
which is periodic in $\theta$, one has
\begin{equation}\label{eq:local-averaging-variable-phase}
\int_B
G\!\left(
x,\frac{\phi(t,x)}{\varepsilon^p},h
\right)\,dx
\longrightarrow
\int_B\int_{\mathbb T}
G(x,\theta,h)\,d\theta\,dx
\end{equation}
as $\varepsilon\to0$, uniformly for $h\in\mathcal H$.
\end{lem}

\begin{Proof}
Since $\nabla_x\phi(t,x_0)\neq0$, after relabeling the coordinates and
shrinking $B$ if necessary, we may assume
$\partial_{x_1}\phi(t,x)\neq0$ on $\overline B$. Hence
\[
\Psi_t(x):=
\bigl(\phi(t,x),x_2,\ldots,x_d\bigr)
\]
is a $C^1$ diffeomorphism from $B$ onto its image
$\widetilde B:=\Psi_t(B)$.

\medskip
\noindent
Writing $y=\Psi_t(x)$, with inverse $x=X_t(y)$, and setting
$J_t(y):=|\det DX_t(y)|$, we obtain
\[
\int_B
G\!\left(
x,\frac{\phi(t,x)}{\varepsilon^p},h
\right)\,dx
=
\int_{\widetilde B}
G\!\left(
X_t(y),\frac{y_1}{\varepsilon^p},h
\right)
J_t(y)\,dy.
\]
The function
$\widetilde G(y,\theta,h):=G(X_t(y),\theta,h)J_t(y)$
is continuous and periodic in $\theta$. The standard one-dimensional
periodic averaging property in the variable $y_1$, uniformly for
$h\in\mathcal H$, therefore gives \cref{eq:local-averaging-variable-phase}
after integration in the remaining variables.
\end{Proof}

\subsubsection{Proof of Lemma~\ref{lem:appendix-wavelength-increment}}
\label{subsubsec:appendix-proof-wavelength-lemma}

\begin{Proof}
Fix $x_0\in\Omega_t$. Since
$\theta\mapsto U^1(t,x_0,\theta)$ is nonconstant, there exists
$\alpha_0\in\mathbb R\setminus\{0\}$ such that
\begin{equation}\label{eq:appendix-positive-shift}
\int_{\mathbb T}
\left|
U^1(t,x_0,\theta+\alpha_0)
-
U^1(t,x_0,\theta)
\right|
\,d\theta
>0.
\end{equation}
Since $\nabla_x\phi(t,x_0)\neq0$, we may choose
$h_0\in\mathbb R^d\setminus\{0\}$ such that
\begin{equation}\label{eq:appendix-choice-h0}
\nabla_x\phi(t,x_0)\cdot h_0
=
\alpha_0.
\end{equation}

\medskip
\noindent
Define
\begin{equation}\label{eq:appendix-G}
G(x,\theta,h)
=
\left|
U^1\!\left(
t,x,\theta+\nabla_x\phi(t,x)\cdot h
\right)
-
U^1(t,x,\theta)
\right|.
\end{equation}
By \cref{eq:appendix-positive-shift,eq:appendix-choice-h0},
$\int_{\mathbb T}G(x_0,\theta,h_0)\,d\theta>0$. By continuity, after
shrinking the ball furnished by
Lemma~\ref{lem:local-averaging-variable-phase} if necessary, there exist a
ball $B_0\Subset\Omega_t$ containing $x_0$, a compact neighborhood
$\mathcal H\Subset\mathbb R^d\setminus\{0\}$ of $h_0$ with non-empty
interior, and $c_0>0$ such that
\begin{equation}\label{eq:appendix-positive-average}
\int_{\mathbb T}
G(x,\theta,h)\,d\theta
\geq
c_0,
\qquad
x\in B_0,\quad h\in\mathcal H.
\end{equation}
Choose $B$ with $B_0\Subset B\Subset\Omega_t$. Since $\mathcal H$ is
compact, $x+\varepsilon^ph\in B$ for
$x\in B_0$, $h\in\mathcal H$, and all sufficiently small $\varepsilon$.

\medskip
\noindent
Taylor's formula gives, uniformly for
$x\in B_0$ and $h\in\mathcal H$,
\begin{equation}\label{eq:appendix-phase-increment}
\frac{
\phi(t,x+\varepsilon^p h)-\phi(t,x)
}{
\varepsilon^p
}
=
\nabla_x\phi(t,x)\cdot h
+
O(\varepsilon^p).
\end{equation}
Using the smoothness of $U_1$, we obtain
\begin{align}
&
\varepsilon
U^1\!\left(
t,x+\varepsilon^p h,
\frac{\phi(t,x+\varepsilon^p h)}{\varepsilon^p}
\right)
-
\varepsilon
U^1\!\left(
t,x,
\frac{\phi(t,x)}{\varepsilon^p}
\right)
\nonumber\\
&\qquad =
\varepsilon
\Bigg[
U^1\!\left(
t,x,
\frac{\phi(t,x)}{\varepsilon^p}
+
\nabla_x\phi(t,x)\cdot h
\right)
-
U^1\!\left(
t,x,
\frac{\phi(t,x)}{\varepsilon^p}
\right)
\Bigg]
+
o(\varepsilon),
\label{eq:appendix-principal-increment}
\end{align}
uniformly for $x\in B_0$ and $h\in\mathcal H$.

\medskip
\noindent
For every $k\geq2$, the corresponding increment is
$O(\varepsilon^k)$, uniformly on $B_0\times\mathcal H$. Since $K$ is fixed,
the total contribution of the higher profiles is
$O(\varepsilon^2)=o(\varepsilon)$. Hence
\begin{align}
&
\left|
(v_K^\varepsilon-\underline u)
(t,x+\varepsilon^p h)
-
(v_K^\varepsilon-\underline u)(t,x)
\right|
\nonumber\\
&\qquad =
\varepsilon
G\!\left(
x,\frac{\phi(t,x)}{\varepsilon^p},h
\right)
+
o(\varepsilon)
\label{eq:appendix-increment-G}
\end{align}
in $L^1(B_0)$, uniformly for $h\in\mathcal H$.

\medskip
\noindent
Applying Lemma~\ref{lem:local-averaging-variable-phase} on $B_0$ gives
\begin{equation}\label{eq:appendix-G-averaging}
\int_{B_0}
G\!\left(
x,\frac{\phi(t,x)}{\varepsilon^p},h
\right)\,dx
\longrightarrow
\int_{B_0}\int_{\mathbb T}
G(x,\theta,h)\,d\theta\,dx
\end{equation}
uniformly for $h\in\mathcal H$. By
\cref{eq:appendix-positive-average}, the limit is bounded below by
$c_0|B_0|$, and therefore, for all sufficiently small $\varepsilon$,
\[
\int_{B_0}
G\!\left(
x,\frac{\phi(t,x)}{\varepsilon^p},h
\right)\,dx
\geq
\frac{c_0|B_0|}{2},
\qquad
h\in\mathcal H.
\]
Combining this estimate with \cref{eq:appendix-increment-G} gives
\cref{eq:appendix-increment-lower} and completes the proof.
\end{Proof}

\bibliographystyle{plain}
\bibliography{biblio}
\end{document}